\input amstex
\documentstyle{amsppt}
\nopagenumbers
\nologo
%------------------------------------------------------

%
% Title:    Quick Introduction to Tensor Analysis.
% Author:   Ruslan Sharipov
% Comments: The textbook, typeset by AmSTeX, 47 pages, 
%           amsppt style, 13 figures in EPS files
%
%------------------------------------------------------
%           Replacement for output macro definition
%
\global\voffset=30pt
\catcode`@=11
\redefine\output@{%
  \def\break{\penalty-\@M}\let\par\endgraf
  \ifodd\pageno\global\hoffset=110pt\else\global\hoffset=-20pt\fi
  \shipout\vbox{%
    \ifplain@
      \let\makeheadline\relax \let\makefootline\relax
    \else
      \iffirstpage@ \global\firstpage@false
        \let\rightheadline\frheadline
        \let\leftheadline\flheadline
      \else
        \ifrunheads@ %\let\makefootline\relax
        \else \let\makeheadline\relax
        \fi
      \fi
    \fi
    \makeheadline \pagebody \makefootline}%
  \advancepageno \ifnum\outputpenalty>-\@MM\else\dosupereject\fi
}
%---------------------------------------------------------------
\font\cpr=cmr7
\newcount\xnumber
\footline={\xnumber=\pageno
\divide\xnumber by 7
\multiply\xnumber by -7
\advance\xnumber by\pageno
\ifnum\xnumber>0\hfil\else\vtop{\vskip 0.5cm
\noindent\cpr CopyRight \copyright\ Sharipov R.A.,
2004.}\hfil\fi}
%---------------------------------------------------------------
\def\setfirstpage{\global\firstpage@true}
\catcode`\@=\active
%---------------------------------------------------------------
% Redefinition of some font parameters and proclaim style
\fontdimen3\tenrm=3pt
\fontdimen4\tenrm=0.7pt

%---------------------------------------------------------------
\def\leaderfill{\leaders\hbox to 0.3em{\hss.\hss}\hfill}
%---------------------------------------------------------------
\Monograph
\loadbold
\TagsOnRight
\pagewidth{360pt}
\pageheight{606pt}
%---------------------------------------------------------------
% Here are local macro definitions for this document
\def\negskp{\hskip -2pt}
\def\compos{\,\raise 1pt\hbox{$\sssize\circ$} \,}
\def\idop{\operatorname{\bold{id}}}
\def\Cl{\operatorname{Cl}}
\def\divr{\operatorname{div}}
\def\rot{\operatorname{rot}}
\def\grad{\operatorname{grad}}
\def\msum{\operatornamewithlimits{\sum^3\!{\ssize\ldots}\!\sum^3}}
\def\bluesum{\operatornamewithlimits{\blue{\sum}}}
\def\blue#1{#1}
\def\red#1{#1}
\accentedsymbol\TTbe{\Tilde{\Tilde{\bold e}}}
\accentedsymbol\TTS{\Tilde{\Tilde{S}}}
\accentedsymbol\TTT{\Tilde{\Tilde{T}}}
%--------------------------------------
\document
% This is cover page
\vbox to\vsize{
\hphantom{a}
\vskip 680pt
\hbox to 0pt{\kern -10pt\includegraphics{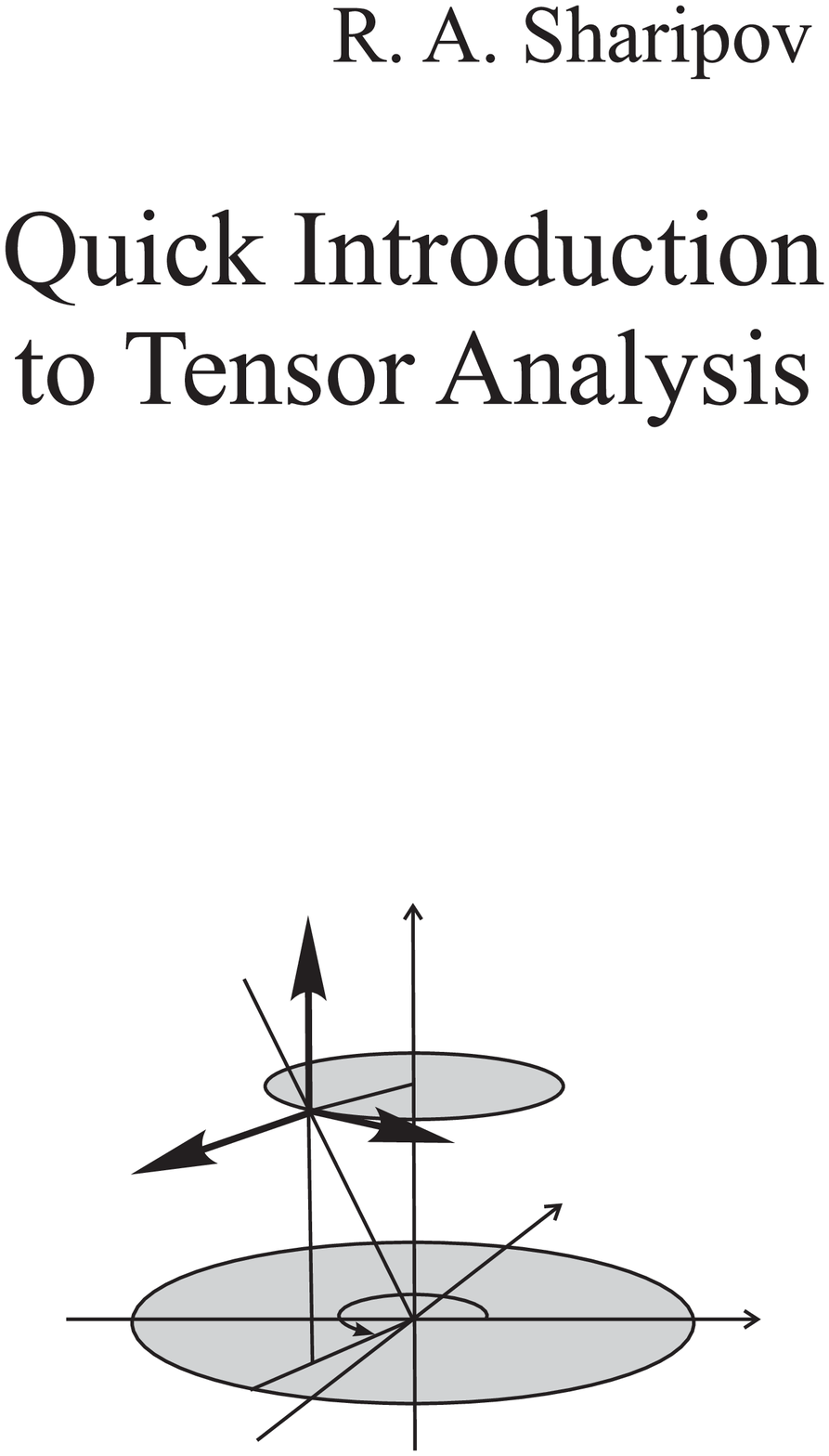}\hss}
\vskip -200pt
\vfill}
%--------------------------------------
% This is second page
\vbox to\vsize{
MSC 97U20\par
\medskip
R. A. Sharipov. {\bf Quick Introduction to Tensor Analysis}:
lecture notes. Freely distributed on-line. Is free for individual
use and educational purposes. Any commercial use without written
consent from the author is prohibited.
\medskip
This book was written as lecture notes for classes that I taught to
undergraduate students majoring in physics in February 2004 during 
my time as a guest instructor at The University of Akron, which was 
supported by Dr. Sergei F. Lyuksyutov's grant from the National 
Research Council under the COBASE program. These 4 classes have been 
taught in the frame of a regular Electromagnetism course as an 
introduction to tensorial methods.\par
     I wrote this book in a "do-it-yourself" style so that I give 
only a draft of tensor theory, which includes formulating definitions 
and theorems and giving basic ideas and formulas. All other work such 
as proving consistence
of definitions, deriving formulas, proving theorems or completing
details to proofs is left to the reader in the form of numerous 
exercises. I hope that this style makes learning the subject really 
quick and more effective for understanding and memorizing.\par
    I am grateful to Department Chair Prof\.~Robert~R.~Mallik for
the opportunity to teach classes and thus to be involved fully in 
the atmosphere of an American university. I am also grateful to 
\roster
\item"" Mr. M.~Boiwka (%
\blue{mboiwka\@hotmail.com})
\item"" Mr. A.~Calabrese (%
\blue{ajc10\@uakron.edu})
\item"" Mr. J.~Comer (%
\blue{funnybef\@lycos.com})
\item"" Mr. A.~Mozinski (%
\blue{arm5\@uakron.edu})
\item"" Mr. M.~J.~Shepard (%
\blue{sheppp2000\@yahoo.com})
\endroster
for attending my classes and reading the manuscript of this book.
I would like to especially acknowledge and thank Mr. Jeff Comer 
for correcting the grammar and wording in it.
\medskip
%\parshape 2 0pt 360pt 50pt 310pt
\noindent {\bf Contacts to author}.
\medskip
\line{\vbox{\hsize=300pt\settabs\+\indent
Office:\ &\cr
\+ Office:\hss &Mathematics Department, Bashkir State University,\cr
\+\hss &32 Frunze street, 450074 Ufa, Russia\cr
\+ Phone:\hss &7-(3472)-23-67-18\cr
\+ Fax:\hss   &7-(3472)-23-67-74\cr
\medskip
\+ Home:\hss &5 Rabochaya street, 450003 Ufa, Russia\cr
\+ Phone:\hss &7-(917)-75-55-786\cr
\+ E-mails: &{\catcode`_=11
\catcode`\_=\active}\blue{R\_\hskip 1pt Sharipov\@ic.bashedu.ru},\cr
\+          &
\blue{r-sharipov\@mail.ru},\cr
\+          &{\catcode`_=11
\catcode`\_=\active}\blue{ra\_\hskip 1pt sharipov\@hotmail.com},\cr
\+ URL: &
\blue{http:/\negskp/www.geocities.com/r-sharipov}\cr
}\hss}
\vfill
\line{CopyRight \copyright\ Sharipov R.A., 2004\hss}}
%--------------------------------------
% This is third page
\topmatter
\title
CONTENTS.
\endtitle
\endtopmatter
\document
\line{CONTENTS.\ \leaderfill\ 3.}
\medskip
\line{CHAPTER~\uppercase\expandafter{\romannumeral 1}.
PRELIMINARY INFORMATION.\ \leaderfill\ 4.}
\medskip
\line{\S~1. Geometrical and physical vectors.\ \leaderfill\ 4.}
\line{\S~2. Bound vectors and free vectors.\ \leaderfill\ 5.}
\line{\S~3. Euclidean space.\ \leaderfill\ 8.}
\line{\S~4. Bases and Cartesian coordinates.\ \leaderfill\ 8.}
\line{\S~5. What if we need to change a basis\,?\ \leaderfill\ 12.}
\line{\S~6. What happens to vectors when we change the basis\,?\ \leaderfill\ 15.}
\line{\S~7. What is the novelty about vectors that we learned knowing\hss}
\line{\qquad transformation formula for their coordinates\,?\ \leaderfill\ 17.}
\setfirstpage
\medskip
\line{CHAPTER~\uppercase\expandafter{\romannumeral 2}.
TENSORS IN CARTESIAN COORDINATES.\ \leaderfill\ 18.}
\medskip
\line{\S~8. Covectors.\ \leaderfill\ 18.}
\line{\S~9. Scalar product of vector and covector.\ \leaderfill\ 19.}
\line{\S~10. Linear operators.\ \leaderfill\ 20.}
\line{\S~11. Bilinear and quadratic forms.\ \leaderfill\ 23.}
\line{\S~12. General definition of tensors.\ \leaderfill\ 25.}
\line{\S~13. Dot product and metric tensor.\ \leaderfill\ 26.}
\line{\S~14. Multiplication by numbers and addition.\ \leaderfill\ 27.}
\line{\S~15. Tensor product.\ \leaderfill\ 28.}
\line{\S~16. Contraction.\ \leaderfill\ 28.}
\line{\S~17. Raising and lowering indices.\ \leaderfill\ 29.}
\line{\S~18. Some special tensors and some useful formulas.\ \leaderfill\ 29.}
\medskip
\line{CHAPTER~\uppercase\expandafter{\romannumeral 3}.
TENSOR FIELDS. DIFFERENTIATION OF TENSORS.\ \leaderfill\ 31.}
\medskip
\line{\S~19. Tensor fields in Cartesian coordinates.\ \leaderfill\ 31.}
\line{\S~20. Change of Cartesian coordinate system.\ \leaderfill\ 32.}
\line{\S~21. Differentiation of tensor fields.\ \leaderfill\ 34.}
\line{\S~22. Gradient, divergency, and rotor. Laplace and d'Alambert 
operators.\ \leaderfill\ 35.}
\medskip
\line{CHAPTER~\uppercase\expandafter{\romannumeral 4}.
TENSOR FIELDS IN CURVILINEAR COORDINATES.\ \leaderfill\ 38.}
\medskip
\line{\S~23. General idea of curvilinear coordinates.\ \leaderfill\ 38.}
\line{\S~24. Auxiliary Cartesian coordinate system.\ \leaderfill\ 38.}
\line{\S~25. Coordinate lines and the coordinate grid.\ \leaderfill\ 39.}
\line{\S~26. Moving frame of curvilinear coordinates.\ \leaderfill\ 41.}
\line{\S~27. Dynamics of moving frame.\ \leaderfill\ 42.}
\line{\S~28. Formula for Christoffel symbols.\ \leaderfill\ 42.}
\line{\S~29. Tensor fields in curvilinear coordinates.\ \leaderfill\ 43.}
\line{\S~30. Differentiation of tensor fields in curvilinear 
coordinates.\ \leaderfill\ 44.}
\line{\S~31. Concordance of metric and connection.\ \leaderfill\ 46.}
\medskip
\line{REFERENCES.\ \leaderfill\ 47.}
\medskip
\newpage
%--------------------------------------
\setfirstpage
\topmatter
\title\chapter{1}
PRELIMINARY INFORMATION.
\endtitle
\endtopmatter
\document
\head
\S~1. Geometrical and physical vectors.
\endhead
\leftheadtext{CHAPTER~\uppercase\expandafter{\romannumeral 1}.
PRELIMINARY INFORMATION.}
    Vector is usually understood as a segment of straight line
equipped with an arrow. Simplest example is displacement vector
$\bold a$. Say its length is $4\,\text{cm}$, i.\,e\.
$$
|\bold a|=4\,\text{cm}.
$$
\setfirstpage
You can draw it on the paper as shown on Fig\.~1a. Then
it means that point $B$ is $4\,\text{cm}$ apart from the
point $A$ in the direction pointed to by vector
$\bold a$. \vadjust{\vskip 130pt\hbox to 0pt{\kern 0pt
\includegraphics{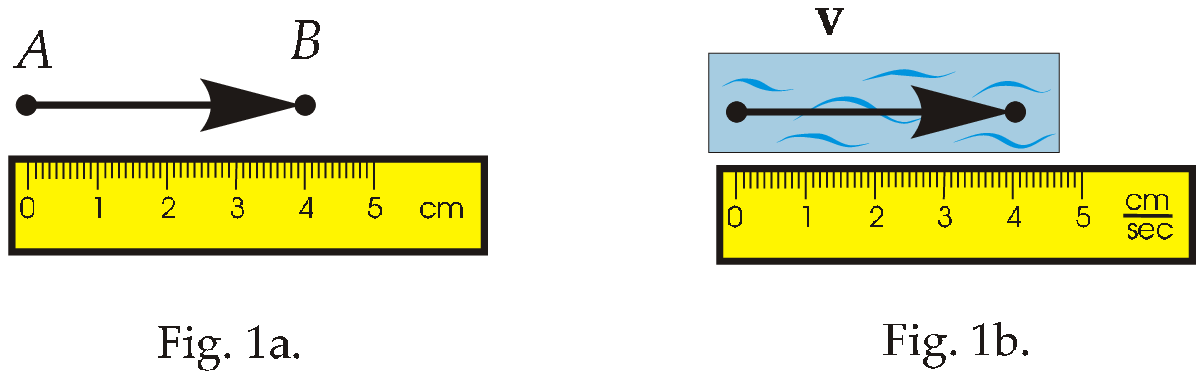}\hss}\vskip -7pt}However,
if you take velocity vector $\bold v$ for a stream in
a brook, you cannot draw it on the paper immediately. You
should first adopt a scaling convention, for example, saying
that $1\,\text{cm}$ on paper represents % corrected by J.C.
$1\,\text{cm/sec}$
(see Fig\.~1b).\par
\proclaim{Conclusion 1.1} Vectors with physical meaning other than
displacement vectors have no unconditional geometric presentation.
Their geometric presentation is conventional; % corrected by J.C.
it depends on the % corrected by J.C.
scaling convention we choose.
\endproclaim
\proclaim{Conclusion 1.2} There are % corrected by J.C.
plenty of physical vectors, which
are not geometrically visible, but can be measured
and then drawn % corrected by J.C.
as geometric vectors.
\endproclaim
    One can consider unit vector $\bold m$. Its length is equal to
unity not $1\,\text{cm}$, not $1\,\text{km}$, not $1\,\text{inch}$,
and not $1\,\text{mile}$, but simply number $1$:
$$
|\bold m|=1.
$$
Like physical vectors, unit vector $\bold m$ cannot be drawn without
adopting some scaling convention. The concept % corrected by J.C.
of a % corrected by J.C.
unit vector is a % corrected by J.C.
very convenient one. \pagebreak By multiplying $\bold m$ to various scalar
quantities, we can produce vectorial quantities of various physical
nature: velocity, acceleration, force, torque, % corrected by J.C.
etc.
\proclaim{Conclusion 1.3} Along with geometrical and physical vectors
one can imagine vectors whose length is a number with
no unit of measure. % corrected by J.C.
\endproclaim
\head
\S~2. Bound vectors and free vectors.
\endhead
    All displacement vectors are bound ones. They are bound to
those points whose displacement they represent. Free vectors are
usually those representing global physical parameters, e\.\,g\.
vector of angular velocity $\boldsymbol\omega$ for Earth rotation
about its axis. \vadjust{\vskip 420pt\hbox to 0pt{\kern 0pt
\includegraphics{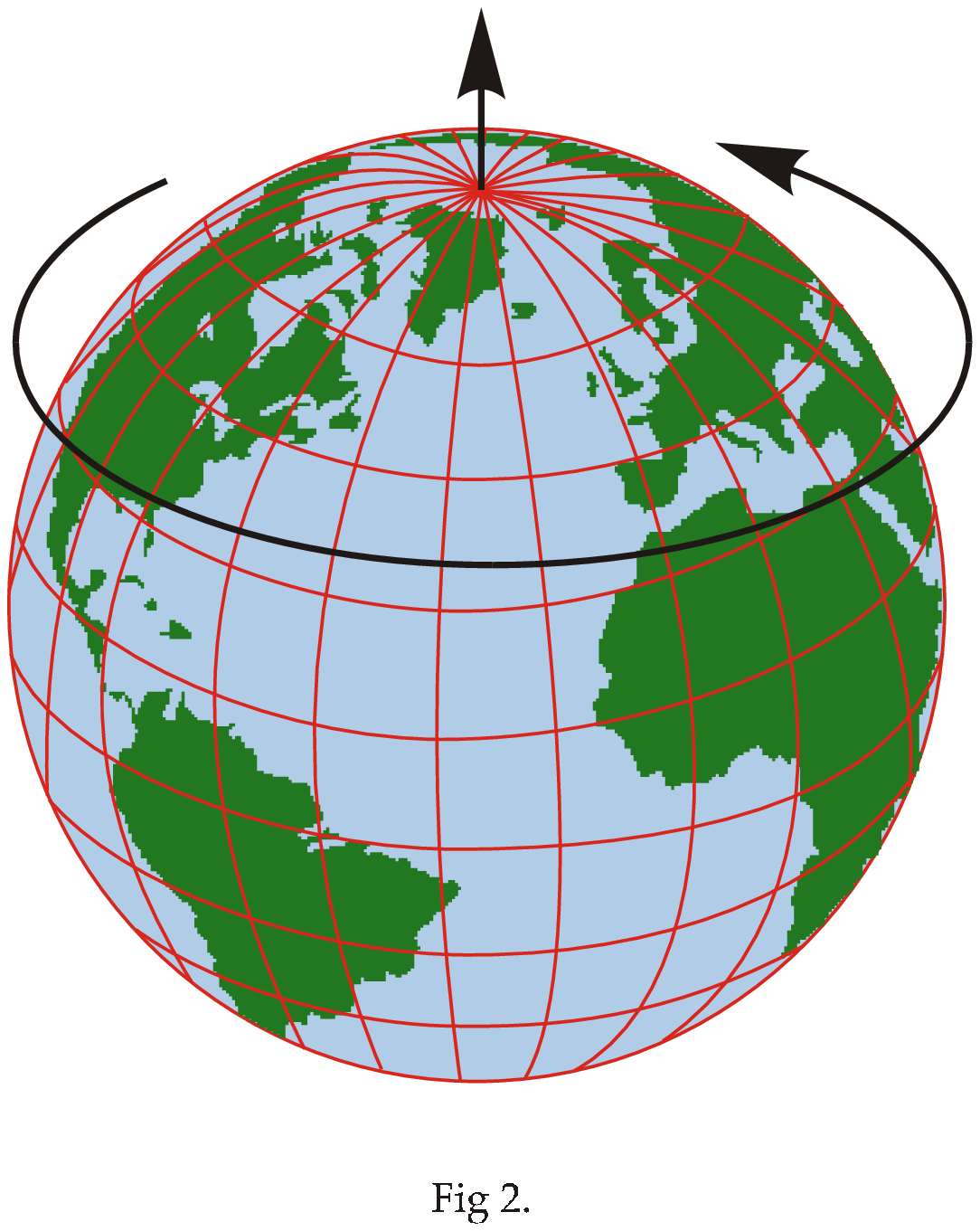}\hss}\vskip -7pt}This vector
produces the % corrected by J.C.
Coriolis force affecting water streams in small rivers and
in oceans around the world. Though it is usually drawn attached to
the % corrected by J.C.
North pole, we can translate this vector to any point along any path
provided we keep its length and \pagebreak direction unchanged.\par
\parshape 3 0pt 360pt 0pt 360pt 180pt 180pt
    The % corrected by J.C.
next example illustrates the concept of a % corrected by J.C.
{\it vector field}.
Consider the water flow in a % corrected by J.C.
river at some fixed instant
of time $t$. \vadjust{\vskip 220pt\hbox to 0pt{\kern 0pt
\includegraphics{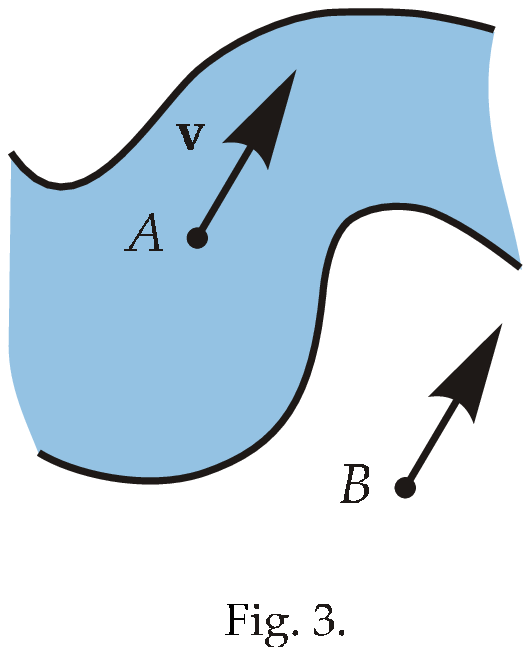}\hss}\vskip -220pt}For each
point $P$ in the water the % corrected by J.C.
velocity of the % corrected by J.C.
water jet passing
through this point is defined. Thus we have a function
$$
\hskip -2em
\bold v=\bold v(t,P).
\tag2.1
$$
Its first argument is time variable $t$. The % corrected by J.C.
second argument of function 
\thetag{2.1} is not numeric. It is geometric object % corrected by J.C.
--- a point. Values of a % corrected by J.C.
function \thetag{2.1} are also not numeric: % corrected by J.C.
they are vectors.
\definition{Definition \!2.1}\parshape 1 180pt 180pt\!A % corrected by J.C.
vector-valued function with point argument is called vector field.
If it has an % corrected by J.C.
additional argument $t$, it is called a % corrected by J.C.
time-dependent vector field.
\enddefinition
\parshape 4 180pt 180pt 180pt 180pt 180pt 180pt 0pt 360pt
    Let $\bold v$ be the value of function \thetag{2.1} at the point
$A$ in a % corrected by J.C.
river. Then vector $\bold v$ is a % corrected by J.C.
bound vector. It represents the % corrected by J.C.
velocity of the % corrected by J.C.
water jet at the point $A$. Hence, % corrected by J.C.
it is bound to % corrected by J.C.
point $A$. Certainly, one can translate it to the point $B$ on the
bank of the river % corrected by J.C.
(see Fig. 3). But there it loses % corrected by J.C.
its original purpose, % corrected by J.C.
which is to mark the water velocity at the point $A$.
\proclaim{Conclusion 2.1} There exist functions with non-numeric
arguments and non-nume\-ric values.
\endproclaim
\proclaim{Exercise 2.1} What is a % corrected by J.C.
scalar field\,? Suggest
an appropriate definition by analogy with definition~2.1.
\endproclaim
\proclaim{Exercise 2.2 {\rm(for deep thinking)}} Let $y=f(x)$
be a function with a % corrected by J.C.
non-numeric argument. Can it be continuous ?
Can it be differentiable ? In general, answer is negative.
However, in some cases one can extend the definition of
continuity and the definition of derivatives in a way applicable
to some functions with non-numeric arguments. Suggest your
version of such a % corrected by J.C.
generalization. If no versions, remember this
problem and return to it later when you gain more experience.
\endproclaim
    Let $A$ be some fixed point (on the ground, under the ground,
in the sky, or in outer space, wherever). % corrected by J.C.
Consider all
vectors of some physical nature bound to this point (say all force
vectors). They constitute an infinite set. Let's denote it
$V_{\sssize A}$. We can perform certain algebraic operations over
the vectors from $V_{\sssize A}$:
\roster
\item we can add any two of them;
\item we can multiply any one of them by any real number
$\alpha\in\Bbb R$;
\endroster
These operations are called linear operations and $V_{\sssize A}$
is called a % corrected by J.C.
linear vector space.
\proclaim{Exercise 2.3} Remember
the parallelogram method % corrected by J.C.
for adding two vectors (draw picture). Remember how vectors are
multiplied by a % corrected by J.C.
real number $\alpha$. Consider three cases:
$\alpha>0$, $\alpha<0$, and $\alpha=0$. Remember what the zero
vector is. % corrected by J.C.
How it is represented geometrically ?
\endproclaim
\proclaim{Exercise 2.4}Do you remember the % corrected by J.C.
exact mathematical definition of a % corrected by J.C.
linear vector space\,? If yes, write it. If no,
visit Web page of Jim Hefferon
\medskip

\centerline{
\rm http:/\negskp/joshua.smcvt.edu/linearalgebra/
}

\medskip
\noindent and download his book \cite{1}.
Keep this book for further references. If you find
it useful, % corrected by J.C.
you can acknowledge the author by sending him e-mail:

{\rm jim\@joshua.smcvt.edu}
.
\endproclaim
\proclaim{Conclusion 2.2} Thus, each point $A$ of our geometric space
is not so simple, even if it is a point in a % corrected by J.C.
vacuum. It can be equipped
with linear vector spaces of various natures % corrected by J.C.
(such as a % corrected by J.C.
space of force vectors in the above example). This idea,
where % corrected by J.C.
each point of vacuum space is
treated as a container for various physical fields, is popular in modern
physics. Mathematically it is realized in the concept of bundles: vector
bundles, tensor bundles, etc.
\endproclaim
\parshape 23 0pt 360pt 0pt 360pt 180pt 180pt 180pt 180pt 180pt 180pt
180pt 180pt 180pt 180pt 180pt 180pt 180pt 180pt  180pt 180pt
180pt 180pt 180pt 180pt 180pt 180pt 180pt 180pt 180pt 180pt
180pt 180pt 180pt 180pt 180pt 180pt 180pt 180pt 180pt 180pt
180pt 180pt 180pt 180pt 0pt 360pt 
     Free vectors, taken as they are, do not form a % corrected by J.C.
linear vector space.
Let's denote by $V$ the set of all free vectors. Then $V$ is union of
vector spaces $V_{\sssize A}$ associated with all points $A$ in space:
$$
\hskip -2em
V=\bigcup_{A\in E}V_{\sssize A}.
\tag2.2
$$
The free vectors forming this set \thetag{2.2} are too
numerous: \vadjust{\vskip 265pt\hbox to 0pt{\kern 0pt
\includegraphics{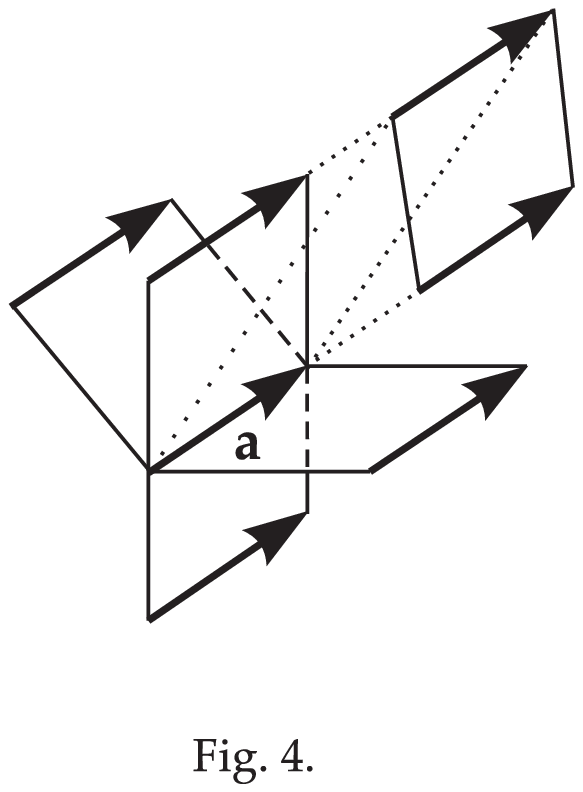}\hss}
\vskip -265pt}we should work to make them 
confine the definition of a linear vector space. Indeed,
if we have a vector $\bold a$ and if it is a free vector,
we can replicate it
by parallel translations that produce % corrected by J.C.
infinitely many copies
of it (see Fig\.~4). All these clones of vector $\bold a$
form a class, the class of vector $\bold a$. Let's denote it
as $\Cl(\bold a)$. Vector $\bold a$ is a representative of its
class. However, we can choose any other vector of this class
as a representative, say it can be vector $\tilde{\bold a}$.
Then we have
$$
\Cl(\bold a)=\Cl(\tilde{\bold a}).
$$
Let's treat $\Cl(\bold a)$ as a whole unit, as one indivisible
object. Then consider the set of all such objects. This set is
called a % corrected by J.C.
factor-set, or quotient set. It is denoted as
$V/\!\sim\,$. This quotient set $V/\!\sim\,$ satisfies the
definition of linear vector space. For the sake of simplicity
further we shall denote it by the same letter $V$ as original
set \thetag{2.2}, from which it is produced by the operation
of factorization.
\proclaim{Exercise 2.5} Have you heard % corrected by J.C.
about binary relations, quotient sets, quotient groups,
quotient rings
and so on\,? If yes, try to remember strict mathematical
definitions for them. If not, then have a look to the 
references \cite{2}, \cite{3}, \cite{4}.
Certainly, you shouldn't read all of these references, but
remember that they are freely available on demand.
\endproclaim
\head
3. Euclidean space.
\endhead
\parshape 3 0pt 360pt 0pt 360pt 175pt 185pt
    What is our geometric space\,? Is it a linear vector space\,?
By no means.  It is formed by points, not by vectors. Properties
of our space were first systematically
described by Euclid, % corrected by J.C.
the Greek mathematician of antiquity.
\vadjust{\vskip 257pt\hbox to 0pt{\kern 0pt \includegraphics{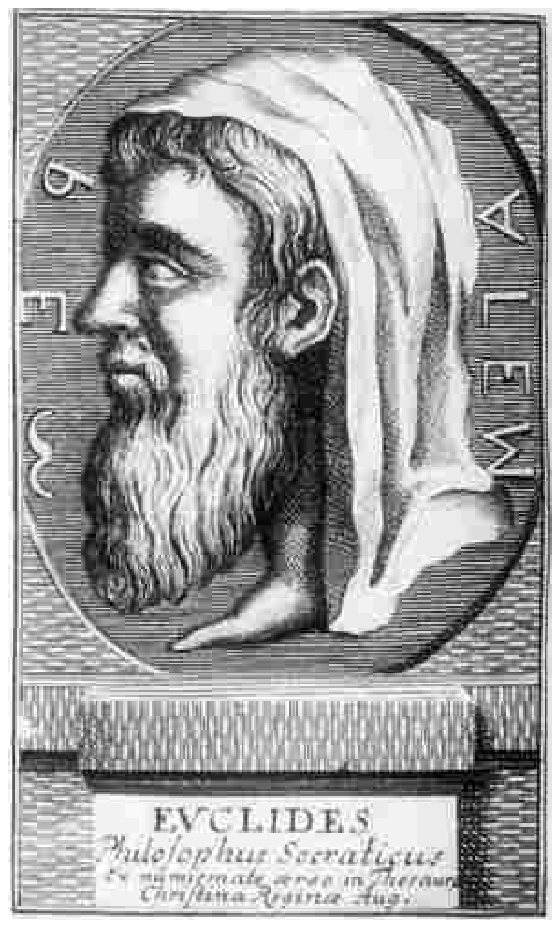}\hss}\vskip -257pt}Therefore, % corrected by J.C.
it is called Euclidean space and denoted by $E$. Euclid
suggested 5 axioms (5 postulates) to describe $E$. However, his
statements were not satisfactorily strict from a % corrected by J.C.
modern point of view.
Currently $E$ is described by 20 axioms. In memory of Euclid
they are subdivided into 5 groups:
\roster
\item\parshape 1 175pt 185pt axioms of incidence;
\item\parshape 1 175pt 185pt axioms of order;
\item\parshape 1 175pt 185pt axioms of congruence;
\item\parshape 1 175pt 185pt axioms of continuity;
\item\parshape 1 175pt 185pt axiom of parallels.
\endroster
\parshape 1 175pt 185pt
20-th axiom, which is also known as 5-th postulate, is most famous.
\proclaim{Exercise 3.1}\!\parshape 1 175pt 185pt Visit the following
{\catcode`~=11
\catcode`\~=\active}\blue{Non-Euclidean Geometry}
web-site and read a few words about non-Euclidean geometry and the
role of Euclid's 5-th postulate in its discovery.
\endproclaim
    Usually nobody remembers % corrected by J.C.
all 20 of these axioms by heart, % corrected by J.C.
even me, though I wrote a textbook on the % corrected by J.C.
foundations of Euclidean geometry in 1998.
Furthermore, % corrected by J.C.
dealing with the Euclidean space $E$, we shall rely only on common
sense and on our geometric intuition.\par
\head
\S~4. Bases and Cartesian coordinates.
\endhead
    Thus, % corrected by J.C.
$E$ is composed by points. Let's choose one of them, denote
it by $O$ and consider the % corrected by J.C.
vector space $V_{\sssize O}$ composed
by displacement vectors. Then each point $B\in E$ can be uniquely
identified with the % corrected by J.C.
displacement vector $\bold r_{\sssize B}=
\overrightarrow{OB}$. It is called the % corrected by J.C.
radius-vector of the point $B$,
while $O$ is called origin. Passing from points to their radius-vectors
we identify $E$ with the % corrected by J.C.
linear vector space $V_{\sssize O}$. Then, passing
from the vectors % corrected by J.C.
to their classes, we can identify $V$ with the space of
free vectors. This identification is a % corrected by J.C.
convenient tool in studying $E$ without referring to
Euclidean % corrected by J.C.
axioms. However, we should remember that such identification
is not unique: it depends on our choice of the
point $O$ for the origin. % corrected by J.C.
\definition{Definition 4.1} We say that three vectors $\bold e_1$,
$\bold e_2$, $\bold e_3$ form a non-coplanar triple of vectors if they
cannot be laid onto the plane by parallel translations.
\enddefinition
These three vectors can be bound to some point $O$ common
to % corrected by J.C.
all of them, or they can be bound to different points in
the space; % corrected by J.C.
it makes no difference. % corrected by J.C.
They also can be treated as free vectors without any
definite binding point.
\definition{Definition 4.2} Any non-coplanar ordered triple of
vectors $\bold e_1$, $\bold e_2$, $\bold e_3$ is called a basis in
our geometric space $E$.
\enddefinition
\proclaim{Exercise 4.1} Formulate the definitions of bases on a plane
and on a % corrected by J.C.
straight line by analogy with definition~4.2.
\endproclaim
\parshape 3 0pt 360pt 0pt 360pt 180pt 180pt
     Below we distinguish three types of bases: orthonormal
basis (ONB), orthogonal basis (OB), and skew-angular basis (SAB).
\vadjust{\vskip 235pt\hbox to 0pt{\kern 0pt\includegraphics{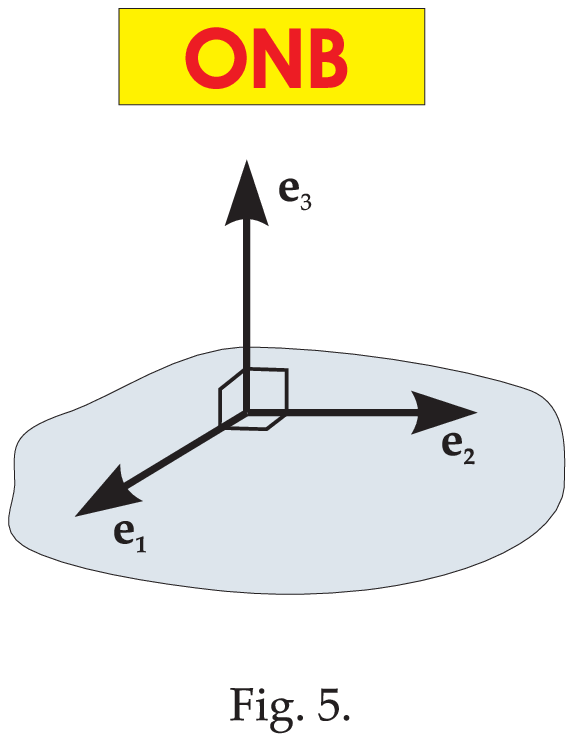}\hss}\vskip -235pt}Orthonormal basis is
formed by three mutually perpendicular unit
vectors: 
$$
\align
\hskip -2em
\aligned
\bold e_1\perp\bold e_2,\\
\bold e_2\perp\bold e_3,\\
\bold e_3\perp\bold e_1,
\endaligned
\tag4.1\\
\vspace{2ex}
\hskip -2em
\aligned
|\bold e_1|=1,\\
|\bold e_2|=1,\\
|\bold e_3|=1.
\endaligned
\tag4.2
\endalign
$$\par
\parshape 1 180pt 180pt
For orthogonal basis, the % corrected by J.C.
three conditions \thetag{4.1} are
fulfilled, but lengths of basis vectors are not specified.\par
\parshape 1 180pt 180pt
    And skew-angular basis is the % corrected by J.C.
most general case. For this basis
neither angles % corrected by J.C.
nor lengths are specified. As we shall see below,
due to its asymmetry SAB can reveal a lot of features
that % corrected by J.C.
are hidden in symmetric ONB.\par
    Let's choose some basis $\bold e_1$, $\bold e_2$, $\bold e_3$
in $E$. In the % corrected by J.C.
general case this is a skew-angular basis. Assume that
vectors $\bold e_1$, $\bold e_2$, $\bold e_3$ are
bound to a % corrected by J.C.
common point $O$ as shown on Fig\.~6 below. Otherwise they can
be brought to this position by means of parallel translations.
Let $\bold a$ be some arbitrary vector. This vector also can be
translated to the point $O$. As a result we have four vectors
$\bold e_1$, $\bold e_2$, $\bold e_3$, and $\bold a$ beginning at
the same point $O$. Drawing additional lines and vectors as shown
on Fig\.~6, we get
$$
\hskip -2em
\bold a=\overrightarrow{OD}=\overrightarrow{OA}+
\overrightarrow{OB}+\overrightarrow{OC}.
\tag4.3
$$
Then from the following obvious relationships
$$
\xalignat 3
&\bold e_1=\overrightarrow{OE_1},
&&\bold e_2=\overrightarrow{OE_2},
&&\bold e_3=\overrightarrow{OE_3},\\
&\overrightarrow{OE_1}\parallel\overrightarrow{OA},
&&\overrightarrow{OE_2}\parallel\overrightarrow{OB},
&&\overrightarrow{OE_3}\parallel\overrightarrow{OC}
\endxalignat 
$$
we derive
$$
\xalignat 3
&\hskip -2em
\overrightarrow{OA}=\alpha\,\bold e_1,
&&\overrightarrow{OB}=\beta\,\bold e_2,
&&\overrightarrow{OC}=\gamma\,\bold e_3,
\tag4.4
\endxalignat
$$
where % corrected by J.C.
$\alpha$, $\beta$, $\gamma$ are scalars. Now from
\thetag{4.3} and \thetag{4.4} we obtain
$$
\hskip -2em
\bold a=\alpha\,\bold e_1+\beta\,\bold e_2
+\gamma\,\bold e_3.
\tag4.5
$$
\proclaim{Exercise 4.2} Explain how, for what reasons, and in
what order % corrected by J.C.
additional lines on Fig\.~6 are drawn.
\vadjust{\vskip 290pt\hbox to 0pt{\kern 0pt
\includegraphics{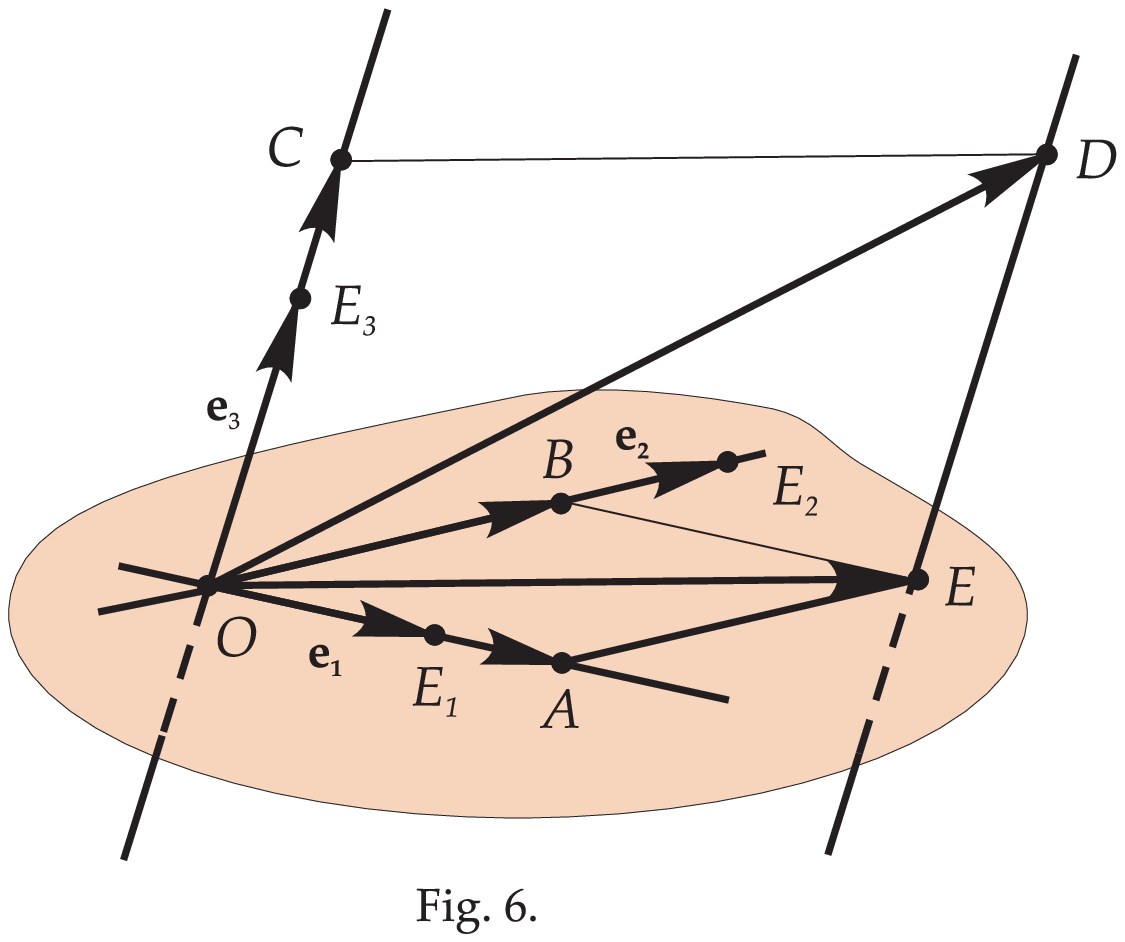}\hss}
\vskip 5pt}
\endproclaim
    Formula \thetag{4.5} is known as the expansion of vector
$\bold a$ in the basis $\bold e_1$, $\bold e_2$, $\bold e_3$,
while $\alpha$, $\beta$, $\gamma$ are coordinates of vector
$\bold a$ in this basis.
\proclaim{Exercise 4.3}Explain why $\alpha$, $\beta$, and
$\gamma$ are uniquely determined by vector $\bold a$.
\endproclaim
\noindent Hint: remember what linear dependence and 
linear independence are. % corrected by J.C.
Give  exact mathematical statements
for these concepts. Apply them to exercise~4.3.\par
    Further we shall write formula \thetag{4.5} as follows
$$
\hskip -2em
\bold a=a^1\,\bold e_1+a^2\,\bold e_2+a^3\,\bold e_3=
\sum^3_{i=1}a^i\,\bold e_i,
\tag4.6
$$
denoting $\alpha=a^1$, $\beta=a^2$, and $\gamma=a^3$.
Don't confuse upper indices in \thetag{4.6} with power
exponentials, $a^1$ here is not $a$, $a^2$ is not $a$
squared, and $a^3$ is not $a$ cubed. Usage of upper
indices and the implicit summation rule were suggested by
Einstein. They are known as Einstein's tensorial
notations.\par
    Once we have chosen the basis $\bold e_1$, $\bold e_2$,
$\bold e_3$ (no matter ONB, OB, or SAB), we can associate
vectors with columns of numbers:
$$
\xalignat 2
&\hskip -2em
\bold a\longleftrightarrow\Vmatrix a^1\\a^2\\a^3\endVmatrix,
&&\bold b\longleftrightarrow\Vmatrix b^1\\b^2\\b^3\endVmatrix.
\tag4.7
\endxalignat
$$
We can % corrected by J.C.
then produce algebraic operations with vectors, reducing them
to arithmetic operations with numbers:
$$
\align
&\bold a+\bold b\longleftrightarrow\Vmatrix a^1\\a^2\\a^3\endVmatrix
+\Vmatrix b^1\\b^2\\b^3\endVmatrix=\Vmatrix a^1+b^1\\a^2+b^2\\a^3+b^3
\endVmatrix,\\
&\alpha\,\bold a\longleftrightarrow\alpha\,\Vmatrix a^1\\a^2\\a^3
\endVmatrix=\Vmatrix\alpha\,a^1\\ \alpha\,a^2\\ \alpha\,a^3\endVmatrix.
\endalign
$$
Columns of numbers framed by matrix delimiters like those in 
\thetag{4.7} are called vector-columns. They form linear vector
spaces.
\proclaim{Exercise 4.4} Remember the % corrected by J.C.
exact mathematical definition for the % corrected by J.C.
real arithmetic vector space $\Bbb R^n$, where $n$ is
a positive integer.
\endproclaim
\definition{Definition 4.1} The % corrected by J.C.
Cartesian coordinate system is a basis
complemented with some fixed point that % corrected by J.C.
is called the % corrected by J.C.
origin.
\enddefinition
     Indeed, if we have an origin $O$, then we can associate each
point $A$ of our space with its radius-vector $\bold r_{\sssize A}=
\overrightarrow{OA}$. Then, having expanded this vector in a basis,
we get three numbers that % corrected by J.C.
are called the Cartesian coordinates of $A$.
Coordinates of a point are also specified % corrected by J.C.
by upper indices since they are coordinates of the % corrected by J.C.
radius-vector for that point. However, unlike coordinates of vectors,
they are usually not written in a column. The reason will be clear when
we consider curvilinear coordinates. So, writing $A(a^1,a^2,a^3)$ is
quite an acceptable notation for the point $A$ with coordinates
$a^1$, $a^2$, and $a^3$.
\vadjust{\vskip 240pt\hbox to 0pt{\kern 0pt
\includegraphics{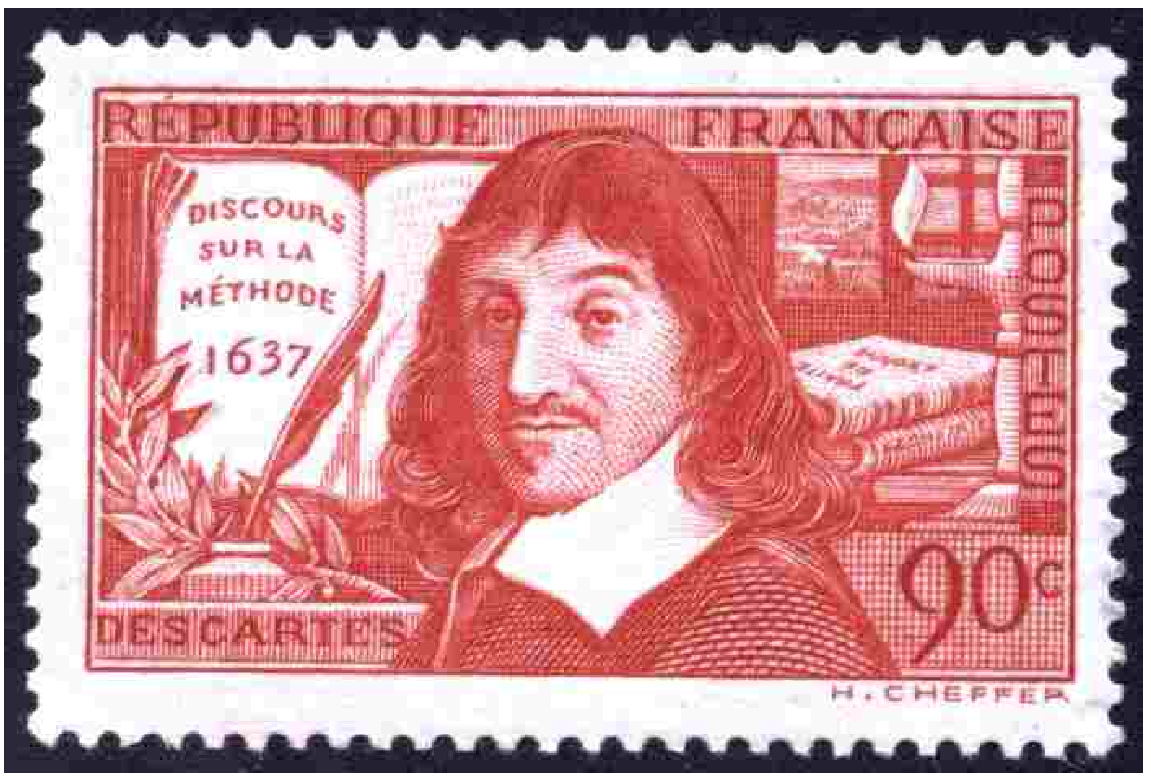}
\hss}\vskip 10pt}\par
    The idea of cpecification of geometric objects by means of
coordinates was first raised by French mathematician and
philosopher Ren\'e Descartes (1596-1650). Cartesian
coordinates are named in memory of him.\par 
\head
\S~5. What if we need to change a basis ?
\endhead
     Why could we need to change a basis\,? There may be
various reasons: % corrected by J.C.
we may dislike initial basis because it is too symmetric like ONB,
or too asymmetric like SAB. Maybe % corrected by J.C.
we are completely satisfied; but the wisdom
is that looking on how something changes we can learn more about this
thing than if we observe it in a % corrected by J.C.
static position. Suppose we have a basis $\bold e_1,\,\bold e_2,\,
\bold e_3$, let's call it {\bf the old basis}, and
suppose we want to change it to a % corrected by J.C.
{\bf new} one $\tilde\bold e_1,\,\tilde
\bold e_2,\,\tilde\bold e_3$. Let's take the % corrected by J.C.
first vector of the % corrected by J.C.
new basis $\bold e_1$. Being isolated from the % corrected by J.C.
other two vectors $\tilde\bold e_2$ and
$\tilde\bold e_3$, it is nothing, but quite an ordinary vector
of space. In this capacity, vector $\tilde\bold e_1$ can be 
expanded in the old basis $\bold e_1,\,\bold e_2,\,\bold e_3$:
$$
\hskip -2em
\tilde\bold e_1=S^1\,\bold e_1+S^2\,\bold e_2+S^3\,\bold e_3=
\sum^3_{j=1}S^j\,\bold e_j.
\tag5.1
$$
Compare \thetag{5.1} and \thetag{4.6}. Then we can take another vector
$\tilde\bold e_2$ and also expand it in the % corrected by J.C.
old basis. But what letter should we % corrected by J.C.
choose for denoting the % corrected by J.C.
coefficients of this expansion\,? We can choose
another letter; say the % corrected by J.C.
letter ``R'':
$$
\hskip -2em
\tilde\bold e_2=R^1\,\bold e_1+R^2\,\bold e_2+R^3\,\bold e_3=
\sum^3_{j=1}R^j\,\bold e_j.
\tag5.2
$$
However, this is not the best decision. Indeed, vectors $\tilde
\bold e_1$ and $\tilde\bold e_2$ differ only in number, while for
their coordinates we use different letters. A % corrected by J.C.
better way is to add
an extra index to $S$ in \thetag{5.1}. This is the % corrected by J.C.
lower index coinciding with the number of the % corrected by J.C.
vector:
$$
\hskip -2em
\tilde\bold e_{\red{1}}
=S^1_{\red{1}}\,\bold e_1+S^2_{\red{1}}\,\bold e_2+S^3_{\red{1}}
\,\bold e_3=\sum^3_{j=1}S^j_{\red{1}}\,\bold e_j
\tag5.3
$$
Color is of no importance; % corrected by J.C.
it is used for highlighting only. Instead of \thetag{5.2},
for the % corrected by J.C.
second vector $\bold e_2$ we write a % corrected by J.C.
formula similar to \thetag{5.3}:
$$
\hskip -2em
\tilde\bold e_{\red{2}}
=S^1_{\red{2}}\,\bold e_1+S^2_{\red{2}}\,\bold e_2+S^3_{\red{2}}
\,\bold e_3=\sum^3_{j=1}S^j_{\red{2}}\,\bold e_j.
\tag5.4
$$
And for third vector as well:
$$
\hskip -2em
\tilde\bold e_{\red{3}}
=S^1_{\red{3}}\,\bold e_1+S^2_{\red{3}}\,\bold e_2+S^3_{\red{3}}
\,\bold e_3=\sum^3_{j=1}S^j_{\red{3}}\,\bold e_j.
\tag5.5
$$
When considered jointly, % corrected by J.C.
formulas \thetag{5.3}, \thetag{5.4},
and \thetag{5.5} are called {\bf transition formulas}.
We use a % corrected by J.C.
left curly bracket to denote their union:
$$
\hskip -2em
\left\{
\aligned
&\tilde\bold e_1=S^1_1\,\bold e_1+S^2_1\,\bold e_2+S^3_1\,\bold e_3,\\
&\tilde\bold e_2=S^1_2\,\bold e_1+S^2_2\,\bold e_2+S^3_2\,\bold e_3,\\
&\tilde\bold e_3=S^1_3\,\bold e_1+S^2_3\,\bold e_2+S^3_3\,\bold e_3.
\endaligned\right.
\tag5.6
$$
We also can write transition formulas
\thetag{5.6} in a % corrected by J.C.
more symbolic form
$$
\hskip -2em
\tilde\bold e_{\red{i}}=\sum^3_{j=1}S^{\blue{j}}_{\red{i}}
\,\bold e_{\blue{j}}.
\tag5.7
$$
Here index $i$ runs over the range of integers % corrected by J.C.
from $1$ to $3$.\par
    Look at index $i$ in formula \thetag{5.7}.
It is a % corrected by J.C.
free index, it can freely take any numeric value from its range:
$1$, $2$, or $3$. Note that $i$ is the % corrected by J.C.
lower index in both sides of formula \thetag{5.7}.
This is a general % corrected by J.C.
rule.
\proclaim{Rule 5.1} In correctly written tensorial formulas free
indices are written on the same level (upper or lower) in both sides
of the equality. Each free index has only one entry in each side of
the equality.
\endproclaim
    Now look at index $j$. It is summation index. It is present
only in right hand side of formula \thetag{5.7}, and it has
exactly two entries (apart from that $j=1$ under the sum symbol):
one in the % corrected by J.C.
upper level and one in the % corrected by J.C.
lower level. This is also general % corrected by J.C.
rule for tensorial formulas.
\proclaim{Rule 5.2} In correctly written tensorial formulas
each summation index should have exactly two entries: one upper
entry and one lower entry.
\endproclaim
    Proposing this rule~5.2, Einstein also suggested
not to write the % corrected by J.C.
summation symbols at all. Formula \thetag{5.7} then would look
like $\tilde\bold e_{\red{i}}=S^{\blue{j}}_{\red{i}}\,
\bold e_{\blue{j}}$ with implicit summation
with respect to the % corrected by J.C.
double index $j$. Many physicists
(especially those % corrected by J.C.
in astrophysics) prefer writing tensorial formulas in this way.
However, I don't like omitting sums.
It breaks the % corrected by J.C.
integrity of notations in science. Newcomers from other branches
of science would have difficulties in understanding formulas with
implicit summation.\par
\proclaim{Exercise 5.1}What happens if $\tilde\bold e_1=
\bold e_1$\,? What are the % corrected by J.C.
numeric values of coefficients $S^1_1$,
$S^2_1$, and $S^3_1$ in formula \thetag{5.3} for this case\,?
\endproclaim
    Returning to transition formulas \thetag{5.6} and \thetag{5.7}
note that coefficients in them are parameterized by two indices
running independently over the range of integer numbers from $1$ to
$3$. In other words, they form a % corrected by J.C.
two-dimensional array that % corrected by J.C.
usually is represented as a table or as a matrix:
$$
\hskip -2em
S=\Vmatrix
S^1_1 & S^1_2 & S^1_3\\
\vspace{1ex}
S^2_1 & S^2_2 & S^2_3\\
\vspace{1ex}
S^3_1 & S^3_2 & S^3_3
\endVmatrix
\tag5.8
$$
Matrix $S$ is called {\bf a % corrected by J.C.
transition matrix} or {\bf direct transition matrix} since
we use it in passing from old basis
to new one. In writing such matrices like $S$ the following
rule applies.
\proclaim{Rule 5.3} For any double indexed array with indices
on the same level (both upper or both lower) the first index
is a row number, while the % corrected by J.C.
second index is a column number. If indices are on different
levels (one upper and one lower), then the % corrected by J.C.
upper index is a row number, while lower one is a column
number.
\endproclaim
    Note that according to this rule~5.3, coefficients of formula
\thetag{5.3}, which are written in line, constitute first column
in matrix \thetag{5.8}. So lines of formula \thetag{5.6} turn
into % corrected by J.C.
columns in matrix \thetag{5.8}. It would be worthwhile % corrected by J.C.
to remember this fact.\par
    If we represent each vector of the % corrected by J.C.
new basis $\tilde\bold e_1,\,\tilde\bold e_2,\,\tilde\bold e_3$
as a column of its coordinates in the % corrected by J.C.
old basis just like it was done for $\bold a$ and $\bold b$ in
formula \thetag{4.7} above
$$
\xalignat 3
&\hskip -2em
\bold e_1\longleftrightarrow\Vmatrix S^1_1\\
\vspace{1ex} S^2_1\\
\vspace{1ex} S^3_1\endVmatrix,
&&\bold e_2\longleftrightarrow\Vmatrix S^1_2\\
\vspace{1ex} S^2_2\\
\vspace{1ex} S^3_2\endVmatrix,
&&\bold e_3\longleftrightarrow\Vmatrix S^1_3\\
\vspace{1ex} S^2_3\\
\vspace{1ex} S^3_3\endVmatrix,
\quad
\tag5.9
\endxalignat
$$
then these columns \thetag{5.9} are exactly the first, the second,
and the third columns in matrix \thetag{5.8}. This is the easiest
way to remember the structure of matrix $S$.\par
\proclaim{Exercise 5.2} What happens if $\tilde\bold e_1=\bold e_1$,
$\tilde\bold e_2=\bold e_2$, and $\tilde\bold e_3=\bold e_3$ ? Find
the % corrected by J.C.
transition matrix for this case. Consider also the following
two cases and write the % corrected by J.C.
transition matrices for each of them:
\roster
\item $\tilde\bold e_1=\bold e_1$, $\tilde\bold e_2=\bold e_3$,
$\tilde\bold e_3=\bold e_2$;
\item $\tilde\bold e_1=\bold e_3$, $\tilde\bold e_2=\bold e_1$,
$\tilde\bold e_3=\bold e_2$.
\endroster
Explain why the next case is impossible:
\roster
\item[3] $\tilde\bold e_1=\bold e_1-\bold e_2$,
$\tilde\bold e_2=\bold e_2-\bold e_3$,
$\tilde\bold e_3=\bold e_3-\bold e_1$.
\endroster
\endproclaim
    Now let's swap bases. This means that we are going to
consider $\tilde\bold e_1,\,\tilde\bold e_2,\,\tilde\bold e_3$
as the % corrected by J.C.
old basis, $\bold e_1,\,\bold e_2,\,
\bold e_3$ as the % corrected by J.C.
new basis, and study the % corrected by J.C.
inverse transition. All of the above stuff applies to
this situation. However, in writing the % corrected by J.C.
transition formulas \thetag{5.6},
let's use another letter for the coefficients. By tradition
here the letter ``T'' is used:
$$
\hskip -2em
\left\{
\aligned
&\bold e_1=T^1_1\,\tilde\bold e_1+T^2_1\,\tilde\bold e_2+T^3_1
\,\tilde\bold e_3,\\
&\bold e_2=T^1_2\,\tilde\bold e_1+T^2_2\,\tilde\bold e_2+T^3_2
\,\tilde\bold e_3,\\
&\bold e_3=T^1_3\,\tilde\bold e_1+T^2_3\,\tilde\bold e_2+T^3_3
\,\tilde\bold e_3.
\endaligned\right.
\tag5.10
$$
Here is the short symbolic version of transition formulas
\thetag{5.10}:
$$
\hskip -2em
\bold e_i=\sum^3_{j=1}T^j_i\,\tilde\bold e_j.
\tag5.11
$$
Denote by $T$ the % corrected by J.C.
transition matrix constructed on the base of
\thetag{5.10} and \thetag{5.11}.
It is called {\bf the % corrected by J.C.
inverse transition matrix} when compared to the % corrected by J.C.
direct transition matrix $S$:
\vskip -1.5ex
$$
\hskip -2em
(\bold e_1,\,\bold e_2,\,\bold e_3)
\vcenter{\hsize 40pt$$\CD @>S>>\\ \vspace{-1.5em} @<<T<\endCD$$}
(\tilde\bold e_1,\,\tilde\bold e_2,\,\tilde\bold e_3).
\tag5.12
$$
\proclaim{Theorem 5.1} The % corrected by J.C.
inverse transition matrix $T$ in
\thetag{5.12} is the % corrected by J.C.
inverse matrix for the direct transition
matrix $S$, i.\,e\. $T=S^{-1}$.
\endproclaim
\proclaim{Exercise 5.3} What is the % corrected by J.C.
inverse matrix\,? Remember the definition. How is the
inverse matrix $A^{-1}$ calculated if $A$ is known\,?
(Don't say that you use a % corrected by J.C.
computer package like
Maple, MathCad, or any other; remember the algorithm
for calculating $A^{-1}$).
\endproclaim
\proclaim{Exercise 5.4} Remember what is the % corrected by J.C.
determinant of a % corrected by J.C.
matrix. How is it % corrected by J.C.
usually calculated\,? Can you calculate $\det(A^{-1})$
if $\det A$ is already known\,?
\endproclaim
\proclaim{Exercise 5.5} What is matrix multiplication\,? Remember
how it is defined. Suppose you have a % corrected by J.C.
rectangular $5\times 3$ matrix $A$ and another rectangular matrix
$B$ which is $4\times 5$. Which of these two products $A\,B$ or
$B\,A$ you can calculate\,?
\endproclaim
\proclaim{Exercise 5.6} Suppose that $A$ and $B$ are two rectangular
matrices, and suppose that $C=A\,B$. Remember the % corrected by J.C.
formula for the components in matrix $C$ if the % corrected by J.C.
components of $A$ and $B$ are known (they are denoted % corrected by J.C.
by $A_{ij}$ and $B_{pq}$). % corrected by J.C.
Rewrite this formula for the case when the % corrected by J.C.
components of $B$ are denoted by $B^{pq}$. Which indices (upper, or
lower, or mixed) you would use for components of $C$ in the last case
(see rules~5.1 and 5.2  of Einstein's tensorial notation).
\endproclaim
\proclaim{Exercise 5.7} Give some examples of matrix multiplication
that % corrected by J.C.
are consistent with Einstein's tensorial notation and those
that % corrected by J.C.
are not (please, do not use examples that % corrected by J.C.
are already considered in exercise~5.6).
\endproclaim
    Let's consider three bases: basis one $\bold e_1,\,\bold e_2,\,
\bold e_3$, basis two $\tilde\bold e_1,\,\tilde\bold e_2,\,\tilde
\bold e_3$, and basis three $\TTbe_1,\,\TTbe_2,\,\TTbe_3$. And let's
consider the % corrected by J.C.
transition matrices relating them:
\vskip -2ex
$$
\hskip -2em
(\bold e_1,\,\bold e_2,\,\bold e_3)
\vcenter{\hsize 40pt$$\CD @>S>>\\ \vspace{-1.5em} @<<T<\endCD$$}
(\tilde\bold e_1,\,\tilde\bold e_2,\,\tilde\bold e_3)
\vcenter{\hsize 40pt$$\CD @>\tilde S>>\\ \vspace{-1.5em}
@<<\tilde T<\endCD$$} (\TTbe_1,\,\TTbe_2,\,\TTbe_3).
\tag5.13
$$
Denote by $\TTS$ and $\TTT$ transition matrices relating basis one
with basis three in \thetag{5.13}:
\vskip -2ex
$$
\hskip -2em
(\bold e_1,\,\bold e_2,\,\bold e_3)
\vcenter{\hsize 40pt$$\CD @>\TTS>>\\ \vspace{-1.5em} @<<\TTT<\endCD$$}
(\TTbe_1,\,\TTbe_2,\,\TTbe_3).
\tag5.14
$$
\proclaim{Exercise 5.8} For matrices $\TTS$ and $\TTT$ in \thetag{5.14}
prove that $\TTS=S\,\tilde S$ and $\TTT=\tilde T\,T$. Apply this result
for proving theorem~5.1.
\endproclaim
\head
\S~6. What happens % corrected by J.C.
to vectors when we change the basis\,?
\endhead
    The answer to this question is very simple. Really nothing\,!
Vectors do not need a basis for their being. But their coordinates, 
they depend on our choice of basis. And they change if we change
the basis. Let's study how % corrected by J.C.
they change. Suppose we have some
vector $\bold x$ expanded in the basis $\bold e_1,\,\bold e_2,\,
\bold e_3$:
$$
\hskip -2em
\bold x=x^1\,\bold e_1+x^2\,\bold e_2+x^3\,\bold e_3
=\sum^3_{i=1}x^i\,\bold e_i.
\tag6.1
$$
Then we keep vector $\bold x$ and change the basis $\bold e_1,
\,\bold e_2,\,\bold e_3$ to % corrected by J.C.
another basis $\tilde\bold e_1,\,
\tilde\bold e_2,\,\tilde\bold e_3$. As we already learned,
this process is described by transition formula \thetag{5.11}:
$$
\bold e_i=\sum^3_{j=1}T^j_i\,\tilde\bold e_j.
$$
Let's substitute this formula into \thetag{6.1} for $\bold e_i$:
$$
\gather
\bold x=\sum^3_{i=1}x^i\!\left(\,\shave{\sum^3_{j=1}T^j_i
\,\bold e_j}\!\right)=\sum^3_{i=1}\sum^3_{j=1}
x^i\,T^j_i\,\tilde\bold e_j=\sum^3_{j=1}\sum^3_{i=1}
x^i\,T^j_i\,\tilde\bold e_j=\\
=\sum^3_{j=1}\!\left(\,\shave{\sum^3_{i=1}T^j_i\,x^i}\!\right)
\,\tilde\bold e_j=\sum^3_{j=1}\tilde x^j\,\tilde\bold e_j
\text{, \ where \ }\tilde x^j=\sum^3_{i=1}T^j_i\,x^i.
\endgather
$$
Thus we have calculated the expansion of vector $\bold x$ in
the % corrected by J.C.
new basis and have derived the % corrected by J.C.
formula relating its new coordinates
to its initial ones:
$$
\hskip -2em
\tilde x^j=\sum^3_{i=1}T^j_i\,x^i.
\tag6.2
$$
This formula is called {\bf a % corrected by J.C.
transformation formula}, or {\bf direct transformation formula}.
Like \thetag{5.7}, it can be written in expanded form:
$$
\hskip -2em
\left\{\aligned
&\tilde x^1=T^1_1\,x^1+T^1_2\,x^2+T^1_3\,x^3,\\
&\tilde x^2=T^2_1\,x^1+T^2_2\,x^2+T^2_3\,x^3,\\
&\tilde x^3=T^3_1\,x^1+T^3_2\,x^2+T^3_3\,x^3.
\endaligned\right.
\tag6.3
$$
And the % corrected by J.C.
transformation formula \thetag{6.2} can be written in
matrix form as well:
$$
\hskip -2em
\Vmatrix\tilde x^1\\
\vspace{1ex}
\tilde x^2\\
\vspace{1ex}
\tilde x^3
\endVmatrix=\Vmatrix
T^1_1 & T^1_2 & T^1_3\\
\vspace{1ex}
T^2_1 & T^2_2 & T^2_3\\
\vspace{1ex}
T^3_1 & T^3_2 & T^3_3
\endVmatrix\,
\Vmatrix x^1\\\vspace{1ex}
x^2\\\vspace{1ex} x^3
\endVmatrix.
\tag6.4
$$
Like \thetag{5.7}, formula \thetag{6.2} can be inverted.
Here is {\bf the % corrected by J.C.
inverse transformation formula} expressing the % corrected by J.C.
initial coordinates of vector $\bold x$ through its new
coordinates:
$$
\hskip -2em
x^j=\sum^3_{i=1}S^j_i\,\tilde x^i.
\tag6.5
$$
\proclaim{Exercise 6.1} By analogy with the above calculations
derive the % corrected by J.C.
inverse transformation formula \thetag{6.5} using formula
\thetag{5.7}.
\endproclaim
\proclaim{Exercise 6.2} By analogy with \thetag{6.3} and
\thetag{6.4} write \thetag{6.5} in expanded form and in
matrix form.
\endproclaim
\proclaim{Exercise 6.3} Derive formula \thetag{6.5} directly
from \thetag{6.2} using the concept of the inverse matrix
$S=T^{-1}$.
\endproclaim
    Note that the % corrected by J.C.
direct transformation formula \thetag{6.2} uses the % corrected by J.C.
inverse transition matrix $T$, and the % corrected by J.C.
inverse transformation formula \thetag{6.5} uses direct transition
matrix $S$. It's funny, but it's really so.
\head
\S~7. What is the novelty about the % corrected by J.C.
vectors that we learned knowing transformation formula for
their coordinates\,?
\endhead
\rightheadtext{\S~7. What is the novelty about the vectors \dots}
    Vectors are too common, % corrected by J.C.
too well-known things for
one to expect that there are some % corrected by J.C.
novelties about them. However, the novelty is that % corrected by J.C.
the method of their treatment can be generalized and then applied to
less customary objects. Suppose, we cannot visually observe vectors
(this is really so for some kinds of them, see section 1), but suppose
we can measure their coordinates in % corrected by J.C.
any basis we choose for this purpose. What then do we know about
vectors\,? And how can we tell them from other (non-vectorial)
objects\,? The answer is in formulas \thetag{6.2} and \thetag{6.5}.
Coordinates of vectors (and only coordinates of vectors) will obey
transformation rules \thetag{6.2} and \thetag{6.5} under a change
of basis. Other objects usually have a % corrected by J.C.
different number % corrected by J.C.
of numeric parameters related to the basis,
and even if they have exactly three coordinates (like vectors have),
their coordinates behave differently under a change of basis.
So transformation formulas \thetag{6.2} and \thetag{6.5} work like
detectors, like a sieve for separating vectors from non-vectors.
What are here % corrected by J.C.
non-vectors, and what kind of geometrical and/or physical
objects of a % corrected by J.C.
non-vectorial nature could exist --- these are questions for
a separate discussion. % corrected by J.C.
Furthermore,  we shall consider only a part
of the set of such % corrected by J.C.
objects, which are called tensors.\par
\newpage
%--------------------------------------
\setfirstpage
\topmatter
\title\chapter{2}
TENSORS IN CARTESIAN COORDINATES.
\endtitle
\endtopmatter
\document
\head
\S~8. Covectors.
\endhead
\leftheadtext{CHAPTER~\uppercase\expandafter{\romannumeral 2}.
TENSORS IN CARTESIAN COORDINATES.}
     In previous 7 sections we learned the following important
fact about vectors: a % corrected by J.C.
vector is a physical object in each basis
of our three-dimensional Euclidean space $E$ represented by three
numbers such that these numbers obey certain transformation
rules when we change the basis. These certain transformation rules
are represented by formulas \thetag{6.2} and \thetag{6.5}.\par
    Now suppose that we have some other physical
object that % corrected by J.C.
is represented by three numbers in each basis, and suppose that
these numbers obey some certain transformation rules when we
change the basis, but these rules are different from \thetag{6.2}
and \thetag{6.5}. Is it possible\,? One can try to
find such an % corrected by J.C.
object in nature. However, in mathematics we have another option.
We can construct such an % corrected by J.C.
object mentally, then study its properties,
and finally look if it is represented somehow in nature.\par
    Let's denote our hypothetical object by $\bold a$, and denote
by $a_1,\,a_2,\,a_3$ that three numbers which represent this object
in the basis $\bold e_1,\,\bold e_2,\,\bold e_3$. By analogy with
vectors we shall call them {\bf coordinates}. But in contrast to
vectors, we intentionally used lower indices when denoting them by
$a_1,\,a_2,\,a_3$. Let's prescribe the following transformation
rules to $a_1,\,a_2,\,a_3$ when we change $\bold e_1,\,\bold e_2,
\,\bold e_3$ to % corrected by J.C.
$\tilde\bold e_1,\,\tilde\bold e_2,\,\tilde\bold e_3$:
$$
\align
&\hskip -2em
\tilde a_j=\sum^3_{i=1}S^i_j\,a_i,
\tag8.1\\
&\hskip -2em
a_j=\sum^3_{i=1}T^i_j\,\tilde a_i.
\tag8.2
\endalign
$$
Here $S$ and $T$ are the same transition matrices as in case of
the % corrected by J.C.
vectors in \thetag{6.2} and \thetag{6.5}. Note that \thetag{8.1}
is sufficient, formula \thetag{8.2} is derived from \thetag{8.1}.
\proclaim{Exercise 8.1} Using the concept of the % corrected by J.C.
inverse matrix $T=S^{-1}$ derive formula \thetag{8.2} from formula
\thetag{8.1}. Compare exercise~8.1 and exercise~6.3.
\endproclaim
\definition{Definition 8.1} A % corrected by J.C.
geometric object $\bold a$
in each basis represented by a % corrected by J.C.
triple of coordinates $a_1,\,a_2,\,a_3$ and such
that its coordinates obey transformation rules \thetag{8.1} and
\thetag{8.2} under a change of basis is called {\bf a % corrected by J.C.
covector}.
\enddefinition
    Looking at the above considerations one can think that we
arbitrarily % corrected by J.C.
chose the transformation formula \thetag{8.1}. However,
this is not so. The choice of the % corrected by J.C.
transformation formula should be self-consistent in the following
sense. Let $\bold e_1,\,\bold e_2,\,\bold e_3$ and $\tilde\bold e_1,
\,\tilde\bold e_2,\,\tilde\bold e_3$
be two bases and let $\TTbe_1,\,\TTbe_2,\,\TTbe_3$ be the third basis
in the space. Let's call them basis one, basis two and basis three for
short. We can pass from basis one to basis three directly, see the
right arrow in \thetag{5.14}. Or we can use basis two as an intermediate
basis, see the right arrows in \thetag{5.13}. In both cases the ultimate
result for the coordinates of a covector in basis three should be the same:
this is the self-consistence requirement. It means that coordinates of
a geometric object should depend on the basis, but not on the way that
they were calculated.
\proclaim{Exercise 8.2} Using \thetag{5.13} and \thetag{5.14},
and relying on the results of exercise~5.8 prove that formulas
\thetag{8.1} and \thetag{8.2} yield a self-consistent way of
defining the\linebreak covector.
\endproclaim
\proclaim{Exercise 8.3} Replace $S$ by $T$ in \thetag{8.1} and
$T$ by $S$ in \thetag{8.2}. Show that the resulting formulas
are not self-consistent.
\endproclaim
   What about the physical reality of covectors\,? Later on we
shall see that covectors do exist in nature. They are the nearest
relatives of vectors. And moreover, we shall see that some well-known
physical objects we thought to be vectors are of covectorial nature
rather than vectorial.
\head
\S~9. Scalar product of vector and covector.
\endhead
    Suppose we have a vector $\bold x$ and a covector $\bold a$. Upon
choosing some basis $\bold e_1,\,\bold e_2,\,\bold e_3$, both of them
have three coordinates: $x^1,\,x^2,\,x^3$ for vector $\bold x$, and
$a_1,\,a_2,\,a_3$ for covector $\bold a$. Let's denote by
$\bigl<\bold a,\,\bold x\bigr>$ the following sum:
$$
\hskip -2em
\bigl<\bold a,\,\bold x\bigr>=\sum^3_{i=1}a_i\,x^i.
\tag9.1
$$
The sum \thetag{9.1} is written in agreement with Einstein's
tensorial notation, see rule~5.2 in section~5 above. It is a
number depending on the vector $\bold x$ and on the covector
$\bold a$. This number is called the scalar product of the
vector $\bold x$ and the covector $\bold a$. We use angular
brackets for this scalar product in order to distinguish it
from the scalar product of two vectors in $E$, which is also
known as the dot product.\par
    Defining the scalar product $\bigl<\bold a,\,\bold x\bigr>$
by means of sum \thetag{9.1} we used the coordinates of vector
$\bold x$ and of covector $\bold a$, which are basis-dependent.
However, the value of sum \thetag{9.1} does not depend on any
basis. Such numeric quantities that do not depend on the choice
of basis are called {\bf scalars} or {\bf true scalars}.
\proclaim{Exercise 9.1} Consider two bases $\bold e_1,\,
\bold e_2,\,\bold e_3$ and $\tilde\bold e_1,\,\tilde\bold e_2,
\,\tilde\bold e_3$, and consider the coordinates of vector
$\bold x$ and covector $\bold a$ in both of them. Relying on
transformation rules \thetag{6.2}, \thetag{6.5}, \thetag{8.1},
and \thetag{8.2} prove the equality
$$
\hskip -2em
\sum^3_{i=1}a_i\,x^i=\sum^3_{i=1}\tilde a_i\,\tilde x^i.
\tag9.2
$$
Thus, you are proving the self-consistence of formula \thetag{9.1}
and showing that the scalar product $\bigl<\bold a,\,\bold x\bigr>$
given by this formula is a true scalar quantity.
\endproclaim
\proclaim{Exercise 9.2} Let $\alpha$ be a real number, let
$\bold a$ and $\bold b$ be two covectors, and let $\bold x$ and
$\bold y$ be two vectors. Prove the following properties of
the scalar product \thetag{9.1}:
\vtop{\hsize=160pt
\roster
\item $\bigl<\bold a+\bold b,\,\bold x\bigr>=
\bigl<\bold a,\,\bold x\bigr>+\bigl<\bold b,\,\bold x\bigr>$;
\item $\bigl<\alpha\,\bold a,\,\bold x\bigr>=
\alpha\,\bigl<\bold a,\,\bold x\bigr>$;
\endroster}
\vtop{\hsize=160pt
\roster
\item[3] $\bigl<\bold a,\,\bold x+\bold y\bigr>=
\bigl<\bold a,\,\bold x\bigr>+\bigl<\bold a,\,\bold y\bigr>$;
\item $\bigl<\bold a,\,\alpha\,\bold x\bigr>=
\alpha\,\bigl<\bold a,\,\bold x\bigr>$.
\endroster}
\endproclaim
\proclaim{Exercise 9.3} Explain why the scalar product
$\bigl<\bold a,\,\bold x\bigr>$ is sometimes called the bilinear
function of vectorial argument $\bold x$ and covectorial argument
$\bold a$. In this capacity, it can be denoted as $f(\bold a,\bold x)$.
Remember our discussion about functions with non-numeric arguments
in section~2.
\endproclaim
\noindent{\bf Important note}. The scalar product \vadjust{\vskip 36pt
\hbox to 0pt{\kern 163pt\includegraphics{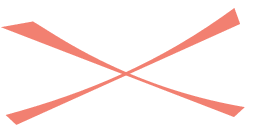}\hss}
\vskip -36pt}$\bigl<\bold a,\,\bold x\bigr>$ is not
symmetric. Moreover, the formula
$$
\bigl<\bold a,\,\bold x\bigr>=\bigl<\bold x,\,\bold a\bigr>
$$
is incorrect in its right hand side since the first argument of scalar
product \thetag{9.1} by definition should be a covector. In a similar
way, the second argument should be a vector. Therefore, we never can
swap them.
\head
\S~10. Linear operators.
\endhead
    In this section we consider more complicated geometric objects.
For the sake of certainty, let's denote one of such objects by
$\bold F$. In each basis $\bold e_1,\,\bold e_2,\,\bold e_3$, it is
represented by a square $3\times 3$ matrix $F^i_j$ of real numbers.
Components of this matrix play the same role as coordinates
in the case of vectors or covectors. Let's prescribe the following
transformation rules to $F^i_j$:
$$
\align
&\hskip -2em
\tilde F^i_j=\sum^3_{p=1}\sum^3_{q=1}T^i_p\,S^q_j\,F^p_q,
\tag10.1\\
&\hskip -2em
F^i_j=\sum^3_{p=1}\sum^3_{q=1}S^i_p\,T^q_j\,\tilde F^p_q.
\tag10.2
\endalign
$$
\proclaim{Exercise 10.1} Using the concept of the inverse matrix
$T=S^{-1}$ prove that formula \thetag{10.2} is derived from
formula \thetag{10.1}.
\endproclaim
If we write matrices $F^i_j$ and $\tilde F^p_q$ according to
the rule~5.3 (see in section 5), then \thetag{10.1} and
\thetag{10.2} can be written as two matrix equalities:
$$
\xalignat 2
&\hskip -2em
\tilde F=T\,F\,S,
&&F=S\,\tilde F\,T.
\tag10.3
\endxalignat
$$
\proclaim{Exercise 10.2} Remember matrix multiplication (we
already considered it in exercises~5.5 and 5.6) and derive
\thetag{10.3} from \thetag{10.1} and \thetag{10.2}.
\endproclaim
\definition{Definition 10.1} A geometric object $\bold F$ in
each basis represented by some square matrix $F^i_j$  and such
that components of its matrix $F^i_j$ obey transformation rules
\thetag{10.1} and \thetag{10.2} under a change of basis is called
{\bf a linear operator}.
\enddefinition
\proclaim{Exercise 10.3} By analogy with exercise~8.2 prove
the self-consistence of the above definition of a linear
operator.
\endproclaim
    Let's take a linear operator $\bold F$ represented by matrix
$F^i_j$ in some basis $\bold e_1,\,\bold e_2,\,\bold e_3$ and take
some vector $\bold x$ with coordinates $x^1,\,x^2,\,x^3$ in the same
basis. Using $F^i_j$ and $x^j$ we can construct the following sum:
$$
\hskip -2em
y^i=\sum^3_{j=1}F^i_j\,x^j.
\tag10.4
$$
Index $i$ in the sum \thetag{10.4} is a free index; it can deliberately
take any one of three values: $i=1$, $i=2$, or $i=3$. For each specific
value of $i$ we get some specific value of the sum \thetag{10.4}.
They are denoted by $y^1,\,y^2,\,y^3$ according to \thetag{10.4}.
Now suppose that we pass to another basis $\tilde\bold e_1,\,
\tilde\bold e_2,\,\tilde\bold e_3$ and do the same things. As a
result we get other three values $\tilde y^1,\,\tilde y^2,\,
\tilde y^3$ given by formula
$$
\hskip -2em
\tilde y^p=\sum^3_{q=1}\tilde F^p_q\,\tilde x^q.
\tag10.5
$$
\proclaim{Exercise 10.4} Relying upon \thetag{10.1} and
\thetag{10.2} prove that the three numbers $y^1,\,y^2,$ $y^3$
and the other three numbers $\tilde y^1,\,\tilde y^2,\,\tilde
y^3$ are related as follows:
$$
\xalignat 2
&\hskip -2em
\tilde y^j=\sum^3_{i=1}T^j_i\,y^i,
&&y^j=\sum^3_{i=1}S^j_i\,\tilde y^i.
\tag10.6
\endxalignat
$$
\endproclaim
\proclaim{Exercise 10.5} Looking at \thetag{10.6} and comparing
it with \thetag{6.2} and \thetag{6.5} find that the $y^1,\,y^2,\,y^3$
and $\tilde y^1,\,\tilde y^2,\,\tilde y^3$ calculated by formulas
\thetag{10.4} and \thetag{10.5} represent the same vector, but in
different bases.
\endproclaim
    Thus formula \thetag{10.4} defines the vectorial object $\bold y$,
while exercise~10.5 assures the correctness of this definition. As a
result we have vector $\bold y$ determined by a linear operator
$\bold F$ and by vector $\bold x$. Therefore, we write
$$
\hskip -2em
\bold y=\bold F(\bold x)
\tag10.7
$$
and say that $\bold y$ is obtained by applying linear operator
$\bold F$ to vector $\bold x$. Some people like to write
\thetag{10.7} without parentheses:
$$
\hskip -2em
\bold y=\bold F\,\bold x.
\tag10.8
$$
Formula \thetag{10.8} is a more algebraistic form of formula
\thetag{10.7}. Here the action of operator $\bold F$ upon vector
$\bold x$ is designated like a kind of multiplication. There is
also a matrix representation of formula \thetag{10.8}, in which
$\bold x$ and $\bold y$ are represented as columns:
$$
\hskip -2em
\Vmatrix
y^1\\\vspace{1ex} y^2\\\vspace{1ex} y^3
\endVmatrix=\Vmatrix
F^1_1 & F^1_2 & F^1_3\\
\vspace{1ex}
F^2_1 & F^2_2 & F^2_3\\
\vspace{1ex}
F^3_1 & F^3_2 & F^3_3
\endVmatrix\,\Vmatrix
x^1\\\vspace{1ex} x^2\\\vspace{1ex} x^3
\endVmatrix.
\tag10.9
$$
\proclaim{Exercise 10.6} Derive \thetag{10.9} from \thetag{10.4}.
\endproclaim
\proclaim{Exercise 10.7} Let $\alpha$ be some real number and
let $\bold x$ and $\bold y$ be two vectors. Prove the following
properties of a linear operator \thetag{10.7}:
\roster
\item $\bold F(\bold x+\bold y)=\bold F(\bold x)+\bold F(\bold y)$,
\item $\bold F(\alpha\,\bold x)=\alpha\,\bold F(\bold x)$.
\endroster
Write these equalities in the more algebraistic style introduced by
\thetag{10.8}. Are they really similar to the properties of
multiplication\,?
\endproclaim
\proclaim{Exercise 10.8} Remember that for the product of two
matrices
$$
\hskip -2em
\det(A\,B)=\det A\,\det B.
\tag10.10
$$
Also remember the formula for $\det(A^{-1})$. Apply these two
formulas to \thetag{10.3} and derive
$$
\hskip -2em
\det F=\det\tilde F.
\tag10.11
$$
\endproclaim
\noindent
Formula \thetag{10.10} means that despite the fact that in
various bases linear operator $\bold F$ is represented by
various matrices, the determinants of all these matrices
are equal to each other. Then we can define the determinant
of linear operator $\bold F$ as the number equal to the
determinant of its matrix in any one arbitrarily chosen
basis $\bold e_1,\,\bold e_2,\,\bold e_3$:
$$
\hskip -2em
\det\bold F=\det F.
\tag10.12
$$
\proclaim{Exercise 10.9 {\rm(for deep thinking)}} Square matrices
have various attributes: eigenvalues, eigenvectors, a characteristic
polynomial, a rank (maybe you remember some others). If we study
these attributes for the matrix of a linear operator, which of them
can be raised one level up and considered as basis-independent
attributes of the linear operator itself\,? Determinant \thetag{10.12}
is an example of such attribute.
\endproclaim
\proclaim{Exercise 10.10} Substitute the unit matrix for $F^i_j$
into \thetag{10.1} and verify that $\tilde F^i_j$ is also a unit
matrix in this case. Interpret this fact.
\endproclaim
\proclaim{Exercise 10.11} Let $\bold x=\bold e_i$ for some
basis $\bold e_1,\,\bold e_2,\,\bold e_3$ in the space.
Substitute this vector $\bold x$ into \thetag{10.7} and by
means of \thetag{10.4} derive the following formula:
$$
\hskip -2em
\bold F(\bold e_i)=\sum^3_{j=1}F^j_i\,\bold e_j.
\tag10.13
$$
Compare \thetag{10.13} and \thetag{5.7}. Discuss the similarities
and differences of these two formulas. The fact is that in some
books the linear operator is determined first, then its matrix is
introduced by formula \thetag{10.13}. Explain why if we know three
vectors $\bold F(\bold e_1)$, $\bold F(\bold e_2)$, and
$\bold F(\bold e_3)$, then we can reconstruct the whole
matrix of operator $\bold F$ by means of formula \thetag{10.13}.
\endproclaim
    Suppose we have two linear operators $\bold F$ and $\bold H$.
We can apply $\bold H$ to vector $\bold x$ and then we can apply
$\bold F$ to vector $\bold H(\bold x)$. As a result we get
$$
\hskip -2em
\bold F\compos\bold H(\bold x)=\bold F(\bold H(\bold x)).
\tag10.14
$$
Here $\bold F\compos\bold H$ is new linear operator introduced by
formula \thetag{10.14}. It is called {\bf a composite operator},
and the small circle sign denotes {\bf composition}.
\proclaim{Exercise 10.12} Find the matrix of composite operator
$\bold F\compos\bold H$ if the matrices for $\bold F$ and $\bold H$
in the basis $\bold e_1,\,\bold e_2,\,\bold e_3$ are known.
\endproclaim
\proclaim{Exercise 10.13} Remember the definition of the identity
map in mathematics
(see 
\blue{\rm on-line Math. Encyclopedia})
and define the identity operator $\idop$. Find the matrix of this
operator.
\endproclaim
\proclaim{Exercise 10.14} Remember the definition of the inverse
map in mathematics and define inverse operator $\bold F^{-1}$ for
linear operator $\bold F$. Find the matrix of this operator if
the matrix of $\bold F$ is known.
\endproclaim
\head
\S~11. Bilinear and quadratic forms.
\endhead
    Vectors, covectors, and linear operators are all examples
of tensors (though we have no definition of tensors yet). Now
we consider another one class of tensorial objects.
For the sake of clarity, let's denote by $a$ one of such objects.
In each basis $\bold e_1,\,\bold e_2,\,\bold e_3$ this object is
represented by some square $3\times 3$ matrix $a_{ij}$ of real
numbers. Under a change of basis these numbers are transformed
as follows:
$$
\align
&\hskip -2em
\tilde a_{ij}=\sum^3_{p=1}\sum^3_{q=1}S^p_i\,S^q_j\,a_{pq},
\tag11.1\\
&\hskip -2em
a_{ij}=\sum^3_{p=1}\sum^3_{q=1}T^p_i\,T^q_j\,\tilde a_{pq}.
\tag11.2
\endalign
$$
Transformation rules \thetag{11.1} and \thetag{11.2} can be
written in matrix form:
$$
\xalignat 2
&\hskip -2em
\tilde a=S^{\sssize\top}\,a\,S,
&&a=T^{\sssize\top}\,\tilde a\,T.
\tag11.3
\endxalignat
$$
Here by $S^{\sssize\top}$ and $T^{\sssize\top}$ we denote
the transposed matrices for $S$ and $T$ respectively.
\proclaim{Exercise 11.1} Derive \thetag{11.2} from \thetag{11.1},
then \thetag{11.3} from \thetag{11.1} and \thetag{11.2}.
\endproclaim
\definition{Definition 11.1} A geometric object $a$ in each basis
represented by some square matrix $a_{ij}$ and such that components
of its matrix $a_{ij}$ obey transformation rules \thetag{11.1} and
\thetag{11.2} under a change of basis is called
{\bf a bilinear form}.
\enddefinition
    Let's consider two arbitrary vectors $\bold x$ and $\bold y$.
We use their coordinates and the components of bilinear form $a$
in order to write the following sum:
$$
\hskip -2em
a(\bold x,\bold y)=\sum^3_{i=1}\sum^3_{j=1}a_{ij}\,x^i\,y^j.
\tag11.4
$$
\proclaim{Exercise 11.2} Prove that the sum in the right hand side
of formula \thetag{11.4} does not depend on the basis, i\.\,e\.
prove the equality
$$
\sum^3_{i=1}\sum^3_{j=1}a_{ij}\,x^i\,y^j=
\sum^3_{p=1}\sum^3_{q=1}\tilde a_{pq}\,\tilde x^p\,\tilde y^q.
$$
\endproclaim
This equality means that $a(\bold x,\bold y)$ is a number determined
by vectors $\bold x$ and $\bold y$ irrespective of the choice of basis.
Hence we can treat \thetag{11.4} as a scalar function of two vectorial
arguments.
\proclaim{Exercise 11.3} Let $\alpha$ be some real number, and let
$\bold x$, $\bold y$, and $\bold z$ be three vectors. Prove the
following properties of function \thetag{11.4}:\newline
\vtop{\hsize=160pt
\roster
\item $a(\bold x+\bold y,\bold z)=
a(\bold x,\bold z)+a(\bold y,\bold z)$;
\item $a(\alpha\,\bold x,\bold y)=
\alpha\,a(\bold x,\bold y)$;
\endroster}
\vtop{\hsize=160pt
\roster
\item[3] $a(\bold x,\bold y+\bold z)=
a(\bold x,\bold y)+a(\bold x,\bold z)$;
\item $a(\bold x,\alpha\,\bold y)=
\alpha\,a(\bold x,\bold y)$.
\endroster}\newline
Due to these properties function \thetag{10.4} is called
a bilinear function or a bilinear form. It is linear with
respect to each of its two arguments.
\endproclaim
   Note that scalar product \thetag{9.1} is also a bilinear
function of its arguments. However, there is a crucial difference
between \thetag{9.1} and \thetag{11.4}. The arguments of scalar
product \thetag{9.1} are of a different nature: the first argument
is a covector, the second is a vector. Therefore, we cannot swap
them. In bilinear form \thetag{11.4} we can swap arguments. As a
result we get another bilinear function
$$
b(\bold x,\bold y)=a(\bold y,\bold x).
\tag11.5
$$
The matrices of $a$ and $b$ are related to each other as follows:
$$
\xalignat 2
&\hskip -2em
b_{ij}=a_{ji},
&&b=a^{\sssize\top}.
\tag11.6
\endxalignat
$$
\definition{Definition 11.2} A bilinear form is called symmetric
if $a(\bold x,\bold y)=a(\bold y,\bold x)$.
\enddefinition
\proclaim{Exercise 11.4} Prove the following identity for
a symmetric bilinear form:
$$
\hskip -2em
a(\bold x,\bold y)=\frac{a(\bold x+\bold y,\bold x+\bold y)
-a(\bold x,\bold x)-a(\bold y,\bold y)}{2}.
\tag11.7
$$
\endproclaim
\definition{Definition 11.3} A quadratic form is a scalar function
of one vectorial argument $f(\bold x)$ produced from some bilinear
function $a(\bold x,\bold y)$ by substituting $\bold y=\bold x$:
$$
\hskip -2em
f(\bold x)=a(\bold x,\bold x).
\tag11.8
$$
\enddefinition
Without a loss of generality a bilinear function $a$ in
\thetag{11.8} can be assumed symmetric. Indeed, if $a$ is not
symmetric, we can produce symmetric bilinear function
$$
\hskip -2em
c(\bold x,\bold y)=\frac{a(\bold x,\bold y)
+a(\bold y,\bold x)}{2},
\tag11.9
$$
and then from \thetag{11.8} due to \thetag{11.9} we derive
$$
\hskip -2em
f(\bold x)=a(\bold x,\bold x)=\frac{a(\bold x,\bold x)
+a(\bold x,\bold x)}{2}=c(\bold x,\bold x).
$$
This equality is the same as \thetag{11.8} with $a$ replaced
by $c$. Thus, each quadratic function $f$ is produced by some
symmetric bilinear function $a$. And conversely, comparing
\thetag{11.8} and \thetag{11.7} we get that $a$ is produced
by $f$:
$$
\hskip -2em
a(\bold x,\bold y)=\frac{f(\bold x+\bold y)
-f(\bold x)-f(\bold y)}{2}.
\tag11.10
$$
Formula \thetag{11.10} is called {\bf the recovery formula}. 
It recovers bilinear function $a$ from quadratic function $f$
produced in \thetag{11.8}. Due to this formula, in referring
to a quadratic form we always imply some symmetric bilinear
form like the geometric tensorial object introduced by
definition~11.1.
\head
\S~12. General definition of tensors.
\endhead
    Vectors, covectors, linear operators, and bilinear forms
are examples of tensors. They are geometric objects that are
represented numerically when some basis in the space is chosen.
This numeric representation is specific to each of them: vectors
and covectors are represented by one-dimensional arrays, linear
operators and quadratic forms are represented by two-dimensional
arrays. Apart from the number of indices, their position does
matter. The coordinates of a vector are numerated by one upper 
index, which is called the contravariant index. The coordinates of
a covector are numerated by one lower index, which is called the 
covariant index. In a matrix of bilinear form we
use two lower indices; therefore bilinear forms are called
{\bf twice-covariant tensors}. Linear operators are tensors
of {\bf mixed type}; their components are numerated by one upper
and one lower index. The number of indices and their positions
determine the transformation rules, i\.\,e\. the way the components
of each particular tensor behave under a change of basis. In the
general case, any tensor is represented by a multidimensional
array with a definite number of upper indices and a definite number
of lower indices. Let's denote these numbers by $r$ and $s$.
Then we have {\bf a tensor of the type $(r,s)$}, or sometimes the
term {\bf valency} is used. A tensor of type $(r,s)$, or of valency
$(r,s)$ is called {\bf an $r$-times contravariant} and
{\bf an $s$-times covariant} tensor. This is terminology; now let's
proceed to the exact definition. It is based on the following general
transformation formulas:
$$
\align
&\hskip -2em
X^{i_1\ldots\,i_r}_{j_1\ldots\,j_s}=\msum\Sb h_1,\,\ldots,\,h_r\\
k_1,\,\ldots,\,k_s\endSb S^{i_1}_{h_1}\ldots\,S^{i_r}_{h_r}T^{k_1}_{j_1}
\ldots\,T^{k_s}_{j_s}\,\tilde X^{h_1\ldots\,h_r}_{k_1\ldots\,k_s},
\tag12.1\\
&\hskip -2em
\tilde X^{i_1\ldots\,i_r}_{j_1\ldots\,j_s}=\msum\Sb h_1,\,\ldots,\,h_r\\
k_1,\,\ldots,\,k_s\endSb T^{i_1}_{h_1}\ldots\,T^{i_r}_{h_r}S^{k_1}_{j_1}
\ldots\,S^{k_s}_{j_s}\,X^{h_1\ldots\,h_r}_{k_1\ldots\,k_s}.
\tag12.2
\endalign
$$
\definition{Definition 12.1} A geometric object $\bold X$ in each
basis represented by $(r+s)$-dimensional array $X^{i_1\ldots\,i_r}_{j_1
\ldots\,j_s}$ of real numbers and such that the components of this
array obey the transformation rules \thetag{12.1} and \thetag{12.2}
under a change of basis is called {\bf tensor} of type $(r,s)$, or
of valency $(r,s)$.
\enddefinition
    Formula \thetag{12.2} is derived from \thetag{12.1}, so it is
sufficient to remember only one of them. Let it be the formula
\thetag{12.1}. Though huge, formula \thetag{12.1} is easy to
remember. One should strictly follow the rules~5.1 and 5.2 from
section~5.\par
    Indices $i_1,\,\ldots,\,i_r$ and $j_1,\,\ldots,\,j_s$
are free indices. In right hand side of the equality \thetag{12.1}
they are distributed in $S$-s and $T$-s, each having only one
entry and each keeping its position, i\.\,e\. upper indices
$i_1,\,\ldots,\,i_r$ remain upper and lower indices $j_1,\,
\ldots,\,j_s$ remain lower in right hand side of the equality
\thetag{12.1}.\par
    Other indices $h_1,\,\ldots,\,h_r$ and $k_1,\,\ldots,\,k_s$
are summation indices; they enter the right hand side of
\thetag{12.1} pairwise: once as an upper index and once as a
lower index, once in $S$-s or $T$-s and once in components of
array $\tilde X^{h_1\ldots\,h_r}_{k_1\ldots\,k_s}$.\par
    When expressing $X^{i_1\ldots\,i_r}_{j_1\ldots\,j_s}$ through
$\tilde X^{h_1\ldots\,h_r}_{k_1\ldots\,k_s}$ each upper index is
served by direct transition matrix $S$ and produces one summation
in \thetag{12.1}:
$$
\hskip -2em
X^{\ldots\,\red{i_{\!\sssize\alpha}}\,\ldots}_{\ldots\,\ldots\,
\ldots}=\sum\ldots\bluesum\limits^3_{h_{\!\sssize\alpha}=1}
\ldots\sum\ldots\,S^{\,\red{i_{\!\sssize\alpha}}}_{\blue{h_{\!
\sssize\alpha}}}\,\ldots\,\tilde X^{\ldots\,\blue{h_{\!\sssize\alpha}}\,
\ldots}_{\ldots\,\ldots\,\ldots}.
\tag12.3
$$
In a similar way, each lower index is served by inverse transition
matrix $T$ and also produces one summation in formula \thetag{12.1}:
$$
\hskip -2em
X^{\ldots\,\ldots\,\ldots}_{\ldots\,\red{j_{\!\sssize\alpha}}\,
\ldots}=\sum\ldots\bluesum\limits^3_{k_{\!\sssize\alpha}=1}
\ldots\sum\ldots\,T^{\blue{k_{\!\sssize\alpha}}}_{\,\red{j_{\!
\sssize\alpha}}}\,\ldots\,\tilde X^{\ldots\,\ldots\,\ldots}_{\ldots\,
\blue{k_{\!\sssize\alpha}}\,\ldots}.
\tag12.4
$$
Formulas \thetag{12.3} and \thetag{12.4} are the same as \thetag{12.1}
and used to highlight how \thetag{12.1} is written. So tensors are
defined. Further we shall consider more examples showing that many
well-known objects undergo the definition~12.1.
\proclaim{Exercise 12.1} Verify that formulas \thetag{6.5},
\thetag{8.2}, \thetag{10.2}, and \thetag{11.2} are special cases
of formula \thetag{12.1}. What are the valencies of vectors, covectors,
linear operators, and bilinear forms when they are considered as
tensors.
\endproclaim
\proclaim{Exercise 12.2} Let $a_{ij}$ be the matrix of some bilinear
form $a$. Let's denote by $b^{ij}$ components of inverse matrix for
$a_{ij}$. Prove that matrix $b^{ij}$ under a change of basis
transforms like matrix of twice-contravariant tensor. Hence it
determines tensor $b$ of valency $(2,0)$. Tensor $b$ is called {\bf
a dual bilinear form} for $a$.
\endproclaim
\head
\S~13. Dot product and metric tensor.
\endhead
    The covectors, linear operators, and bilinear forms that we
considered above were artificially constructed tensors.
However there are some tensors of natural origin. Let's remember
that we live in a space with measure. We can measure distance
between points (hence we can measure length of vectors) and we
can measure angles between two directions in our space.
Therefore for any two vectors $\bold x$ and $\bold y$ we can
define their natural scalar product (or dot product):
$$
\hskip -2em
(\bold x,\,\bold y)=|\bold x|\,|\bold y|\,\cos(\varphi),\
\tag13.1
$$
where $\varphi$ is the angle between vectors $\bold x$ and $\bold y$.
\proclaim{Exercise 13.1} Remember the following properties of the
scalar product \thetag{13.1}:\newline
\vtop{\hsize=160pt\roster
\item $(\bold x+\bold y,\,\bold z)=
(\bold x,\,\bold z)+(\bold y,\,\bold z)$;
\item $(\alpha\,\bold x,\,\bold y)=
\alpha\,(\bold x,\,\bold y)$;
\endroster}
\vtop{\hsize=160pt
\roster
\item[3] $(\bold x,\,\bold y+\bold z)=
(\bold x,\,\bold y)+(\bold x,\,\bold z)$;
\item $(\bold x,\,\alpha\,\bold y)=
\alpha\,(\bold x,\,\bold y)$;
\endroster}
\roster
\item[5] $(\bold x,\,\bold y)=
(\bold y,\,\bold x)$;
\item $(\bold x,\,\bold x)\geqslant 0$ and
$(\bold x,\,\bold x)=0$ implies $\bold x=0$.
\endroster
\vskip 5pt\noindent
These properties are usually considered in courses on
analytic geometry or vector algebra, see

\blue{Vector Lessons on the Web}.
\endproclaim
Note that the first four properties of the scalar product
\thetag{13.1} are quite similar to those for quadratic forms,
see exercise~11.3. This is not an occasional coincidence.
\proclaim{Exercise 13.2} Let's consider two arbitrary vectors
$\bold x$ and $\bold y$ expanded in some basis $\bold e_1,\,
\bold e_2,\,\bold e_3$. This means that we have the following
expressions for them:
$$
\pagebreak
\xalignat 2
&\hskip -2em
\bold x=\sum^3_{i=1}x^i\,\bold e_i,
&&\bold y=\sum^3_{j=1}x^j\,\bold e_j.
\tag13.2
\endxalignat
$$
Substitute \thetag{13.2} into \thetag{13.1} and using properties
\therosteritem{1}--\,\therosteritem{4} listed in exercise~13.1
derive the following formula for the scalar product of $\bold x$ and
$\bold y$:
$$
\hskip -2em
(\bold x,\,\bold y)=\sum^3_{i=1}\sum^3_{j=1}
(\bold e_i,\,\bold e_j)\,x^i\,y^j.
\tag13.3
$$
\endproclaim
\proclaim{Exercise 13.3} Denote $g_{ij}=(\bold e_i,\,\bold e_j)$
and rewrite formula \thetag{13.3} as
$$
\hskip -2em
(\bold x,\,\bold y)=\sum^3_{i=1}\sum^3_{j=1}
g_{ij}\,x^i\,y^j.
\tag13.4
$$
Compare \thetag{13.4} with formula \thetag{11.4}. Consider some
other basis $\tilde\bold e_1,\,\tilde\bold e_2,\,\tilde\bold e_3$,
denote $\tilde g_{pq}=(\tilde\bold e_p,\,\tilde\bold e_q)$ and
by means of transition formulas \thetag{5.7} and \thetag{5.11}
prove that matrices $g_{ij}$ and $\tilde g_{pq}$ are components
of a geometric object obeying transformation rules \thetag{11.1}
and \thetag{11.2} under a change of base. Thus you prove that
the Gram matrix
$$
\hskip -2em
g_{ij}=(\bold e_i,\,\bold e_j)
\tag13.5
$$
defines tensor of type $(0,2)$. This is very important tensor;
it is called {\bf the metric tensor}. It describes not only the
scalar product in form of \thetag{13.4}, but the whole geometry of
our space. Evidences for this fact are below.
\endproclaim
    Matrix \thetag{13.5} is symmetric due to property
\therosteritem{5} in exercise~13.1. Now, comparing \thetag{13.4}
and \thetag{11.4} and keeping in mind the tensorial nature of
matrix \thetag{13.5}, we conclude that the scalar product is a
symmetric bilinear form:
$$
\hskip -2em
(\bold x,\,\bold y)=g(\bold x,\bold y).
\tag13.6
$$
The quadratic form corresponding to \thetag{13.6} is very simple:
$f(\bold x)=g(\bold x,\bold x)=|\bold x|^2$. The inverse matrix
for \thetag{13.5} is denoted by the same symbol $g$ but
with upper indices: $g^{ij}$. It determines a tensor of type
$(2,0)$, this tensor is called {\bf dual metric tensor} (see
exercise~12.2 for more details).
\head
\S~14. Multiplication by numbers and addition.
\endhead
   Tensor operations are used to produce new tensors from those
we already have. The most simple of them are {\bf multiplication by
number} and {\bf addition}. If we have some tensor $\bold X$ of
type $(r,s)$ and a real number $\alpha$, then in some base $\bold
e_1,\,\bold e_2,\,\bold e_3$ we have the array of components
of tensor $X$; let's denote it $X^{i_1\ldots\,i_r}_{j_1\ldots
\,j_s}$. Then by multiplying all the components of this array by
$\alpha$ we get another array
$$
\hskip -2em
Y^{i_1\ldots\,i_r}_{j_1\ldots\,j_s}=\alpha\,
X^{i_1\ldots\,i_r}_{j_1\ldots\,j_s}.
\tag14.1
$$
Choosing another base $\tilde\bold e_1,\,\tilde\bold e_2,\,
\tilde\bold e_3$ and repeating this operation we get
$$
\tilde Y^{i_1\ldots\,i_r}_{j_1\ldots\,j_s}=\alpha\,
\tilde X^{i_1\ldots\,i_r}_{j_1\ldots\,j_s}.
\tag14.2
$$
\proclaim{Exercise 14.1} Prove that arrays $\tilde Y^{i_1\ldots
\,i_r}_{j_1\ldots\,j_s}$ and $Y^{i_1\ldots\,i_r}_{j_1
\ldots\,j_s}$ are related to each other in the same way as
arrays $\tilde X^{i_1\ldots\,i_r}_{j_1\ldots\,j_s}$ and $X^{i_1
\ldots\,i_r}_{j_1\ldots\,j_s}$, i\.\,e\. according to 
transformation formulas \thetag{12.1} and \thetag{12.2}.
In doing this you prove that formula \thetag{14.1} applied in all
bases produces new tensor $\bold Y=\alpha\,\bold X$ from initial
tensor $\bold X$.
\endproclaim
    Formula \thetag{14.1} defines {\bf the multiplication of tensors
by numbers}. In exercise~14.1 you prove its consistence. The next
formula defines {\bf the addition of tensors}:
$$
\hskip -2em
X^{i_1\ldots\,i_r}_{j_1\ldots\,j_s}+
Y^{i_1\ldots\,i_r}_{j_1\ldots\,j_s}=
Z^{i_1\ldots\,i_r}_{j_1\ldots\,j_s}.
\tag14.3
$$
Having two tensors $\bold X$ and $\bold Y$ both of type $(r,s)$
we produce a third tensor $\bold Z$ of the same type $(r,s)$ by
means of formula \thetag{14.3}. It's natural to denote
$\bold Z=\bold X+\bold Y$.
\proclaim{Exercise 14.2} By analogy with exercise~14.1 prove
the consistence of formula \thetag{14.3}.
\endproclaim
\proclaim{Exercise 14.3} What happens if we multiply tensor
$\bold X$ by the number zero and by the number minus one\,?
What would you call the resulting tensors\,?
\endproclaim
\head
\S~15. Tensor product.
\endhead
    The tensor product is defined by a more tricky formula. Suppose
we have tensor $\bold X$ of type $(r,s)$ and tensor $\bold Y$
of type $(p,q)$, then we can write:
$$
\hskip -2em
Z^{i_1\ldots\,i_{r+p}}_{j_1\ldots\,j_{s+q}}=
X^{i_1\ldots\,i_r}_{j_1\ldots\,j_s}\,
Y^{i_{r+1}\ldots\,i_{r+p}}_{j_{s+1}\ldots\,j_{s+q}}.
\tag15.1
$$
Formula \thetag{15.1} produces new tensor $\bold Z$
of the type $(r+p,s+q)$. It is called {\bf the tensor product}
of $\bold X$ and $\bold Y$ and denoted $\bold Z=\bold X
\otimes\bold Y$. Don't mix the tensor product and the cross product.
They are different.
\proclaim{Exercise 15.1} By analogy with exercise~14.1 prove
the consistence of formula \thetag{15.1}.
\endproclaim
\proclaim{Exercise 15.2} Give an example of two tensors
such that $\bold X\otimes\bold Y\neq\bold Y\otimes\bold X$.
\endproclaim
\head
\S~16. Contraction.
\endhead
   As we have seen above, the tensor product increases the number
of indices. Usually the tensor $\bold Z=\bold X\otimes\bold Y$ has
more indices than $\bold X$ and $\bold Y$. Contraction is an
operation that decreases the number of indices. Suppose we have
tensor $\bold X$ of the type $(r+1,s+1)$. Then we can produce tensor
$\bold Z$ of type $(r,s)$ by means of the following formula:
$$
\hskip -2em
Z^{i_1\ldots\,i_r}_{j_1\ldots\,j_s}=\sum^n_{\rho=1}
X^{i_1\ldots\,i_{m-1}\,\red{\rho}\,i_m\ldots\,i_r}_{j_1\ldots
\,j_{k-1}\,\red{\rho}\,j_k\ldots\,j_s}.
\tag16.1
$$
What we do\,? Tensor $\bold X$ has at least one upper index
and at least one lower index. We choose the $m$-th upper index
and replace it by the summation index $\rho$. In the same way we
replace the $k$-th lower index by $\rho$. Other $r$ upper
indices and $s$ lower indices are free. They are numerated
in some convenient way, say as in formula \thetag{16.1}.
Then we perform summation with respect to index $\rho$.
The contraction is over. This operation is called
{\bf a contraction with respect to {\rm $m$-th} upper and
{\rm $k$-th} lower indices}.
Thus, if we have many upper an many lower indices in tensor $\bold
X$, we can perform various types of contractions to this tensor.
\proclaim{Exercise 16.1} Prove the consistence of formula \thetag{16.1}.
\endproclaim
\proclaim{Exercise 16.2} Look at formula \thetag{9.1} and interpret
this formula as the contraction of the tensor product $\bold a\otimes
\bold x$. Find similar interpretations for \thetag{10.4}, \thetag{11.4},
and \thetag{13.4}.
\endproclaim
\head
\S~17. Raising and lowering indices.
\endhead
   Suppose that $\bold X$ is some tensor of type $(r,s)$. Let's
choose its $\alpha$-th lower index: $X^{\ldots\,\ldots\,
\ldots}_{\ldots\,\blue{k}\,\ldots}$. The symbols used for the other
indices are of no importance. Therefore, we denoted them by dots.
Then let's consider the tensor product $\bold Y=g\otimes\bold X$:
$$
\hskip -2em
Y^{\ldots\,p\,q\,\ldots}_{\ldots\,\blue{k}\,\ldots}=
g^{pq}\,X^{\ldots\,\ldots\,\ldots}_{\ldots\,\blue{k}\,\ldots}.
\tag17.1
$$
Here $g$ is the dual metric tensor with the components $g^{pq}$
(see section~13 above). In the next step let's contract \thetag{17.1}
with respect to the pair of indices $k$ and $q$. For this purpose
we replace them both by $s$ and perform the summation:
$$
\hskip -2em
X^{\ldots\,p\,\ldots}_{\ldots\,\ldots\,\ldots}
=\sum^3_{s=1}g^{p\blue{s}}\,X^{\ldots\,\ldots\,\ldots}_{\ldots
\,\blue{s}\,\ldots}.
\tag17.2
$$
This operation \thetag{17.2} is called {\bf the index raising procedure}.
It is invertible. The inverse operation is called {\bf the index lowering
procedure}:
$$
\hskip -2em
X^{\ldots\,\ldots\,\ldots}_{\ldots\,p\,\ldots}
=\sum^3_{s=1}g_{p\blue{s}}\,X^{\ldots\,\blue{s}\,
\ldots}_{\ldots\,\ldots\,\ldots}.
\tag17.3
$$
Like \thetag{17.2}, the index lowering procedure \thetag{17.3} comprises
two tensorial operations: the tensor product and contraction.
\head
\S~18. Some special tensors and some useful formulas.
\endhead
    Kronecker symbol is a well known object. This is a two-dimensional
array representing the unit matrix. It is determined as follows:
$$
\hskip -2em
\delta^i_j=\cases 1 &\text{\ \ for \ }i=j,\\
0 &\text{\ \ for \ }i\neq j.\endcases
\tag18.1
$$
We can determine two other versions of Kronecker symbol:
$$
\hskip -2em
\delta^{ij}=\delta_{ij}=\cases 1 &\text{\ \ for \ }i=j,\\
0 &\text{\ \ for \ }i\neq j.\endcases
\tag18.2
$$
\proclaim{Exercise 18.1} Prove that definition \thetag{18.1}
is invariant under a change of basis, if we interpret the Kronecker
symbol as a tensor. Show that both definitions in \thetag{18.2}
are not invariant under a change of basis.
\endproclaim
\proclaim{Exercise 18.2} Lower index $i$ of tensor \thetag{18.1}
by means of \thetag{17.3}. What tensorial object do you get as
a result of this operation\,?
\endproclaim
\proclaim{Exercise 18.3} Likewise, raise index $J$ in \thetag{18.1}.
\endproclaim
    Another well known object is the Levi-Civita symbol. This is
a three-dimensional array determined by the following formula:
$$
\hskip -2em
\epsilon_{jkq}=\epsilon^{jkq}=
\cases
0,&\vtop{\hsize=4.1cm\noindent\baselineskip 0pt if among 
          $j$, $k$, $q$, there are at least two equal
          numbers;}\\
1,&\vtop{\hsize=4.1cm\noindent\baselineskip 0pt if $(j\,k\,q)$
          is even permutation of numbers $(1\,2\,3)$;}\\
-1,&\vtop{\hsize=4.1cm\noindent\baselineskip 0pt if $(j\,k\,q)$
          is odd permutation of numbers $(1\,2\,3)$.}
\endcases
\tag18.3
$$
The Levi-Civita symbol \thetag{18.3} is not a tensor. However, we
can produce two tensors by means of Levi-Civita symbol. The first
of them
$$
\hskip -2em
\omega_{ijk}=\sqrt{{\det(g_{ij})}}\,\epsilon_{ijk}
\tag18.4
$$
is known as {\bf the volume tensor}. Another one is {\bf the dual
volume tensor}:
$$
\hskip -2em
\omega^{ijk}=\sqrt{\det(g^{ij})\vphantom{g_{ij}}}\,\epsilon^{ijk}.
\tag18.5
$$
Let's take two vectors $\bold x$ and $\bold y$. Then using \thetag{18.4}
we can produce covector $\bold a$:
$$
a_i=\sum^3_{j=1}\sum^3_{k=1}\omega_{ijk}\,x^j\,y^k.
\tag18.6
$$
Then we can apply index raising procedure \thetag{17.2} and
produce vector $\bold a$:
$$
a^r=\sum^3_{i=1}\sum^3_{j=1}\sum^3_{k=1}g^{ri}\,\omega_{ijk}
\,x^j\,y^k.
\tag18.7
$$
Formula \thetag{18.7} is known as formula for the vectorial product
(cross product) in skew-angular basis.
\proclaim{Exercise 18.4} Prove that the vector $\bold a$ with components
\thetag{18.7} coincides with cross product of vectors $\bold x$ and
$\bold y$, i\.\,e\. $\bold a=[\bold x,\,\bold y]$.
\endproclaim
\newpage
%--------------------------------------
\setfirstpage
\topmatter
\title\chapter{3}
TENSOR FIELDS. DIFFERENTIATION OF TENSORS.
\endtitle
\endtopmatter
\document
\head
\S~19. Tensor fields in Cartesian coordinates.
\endhead
\leftheadtext{CHAPTER~\uppercase\expandafter{\romannumeral 3}.
TENSOR FIELDS. DIFFERENTIATION OF TENSORS.}
     The tensors that we defined in section~12 are free tensors.
Indeed, their components are arrays related to bases, while
any basis is a triple of free vectors (not bound to any point).
Hence, the tensors previously considered are also not bound to
any point.
\par
     Now suppose we want to bind our tensor to some point in
space, then another tensor to another point and so on. Doing
so we can fill our space with tensors, one per each point.
In this case we say that we have a tensor field. In order to
mark a point $P$ to which our particular tensor is bound we
shall write $P$ as an argument:
$$
\hskip -2em
\bold X=\bold X(P).
\tag19.1
$$
Usually the valencies of all tensors composing the tensor field
are the same. Let them all be of type $(r,s)$. Then if we choose
some basis $\bold e_1,\,\bold e_2,\,\bold e_3$, we can represent
any tensor of our tensor field as an array $X^{i_1\ldots\,
i_r}_{j_1\ldots\,j_s}$ with $r+s$ indices:
$$
\hskip -2em
X^{i_1\ldots\,i_r}_{j_1\ldots\,j_s}=X^{i_1\ldots\,
i_r}_{j_1\ldots\,j_s}(P).
\tag19.2
$$
Thus, the tensor field \thetag{19.1} is a tensor-valued function
with argument $P$ being a point of three-dimensional Euclidean
space $E$, and \thetag{19.2} is the basis representation for
\thetag{19.1}. For each fixed set of numeric values of indices
$i_1,\,\ldots,\,i_r,\,j_1,\,\ldots,\,j_s$ in \thetag{19.2}, 
we have a numeric function with a point-valued argument. Dealing
with point-valued arguments is not so convenient, for example,
if we want to calculate derivatives. Therefore, we need to replace
$P$ by something numeric. Remember that we have already chosen
a basis. If, in addition, we fix some point $O$ as an origin, then
we get Cartesian coordinate system in space and hence can represent
$P$ by its radius-vector $\bold r_{\sssize\!P}
=\overrightarrow{\vphantom{\tilde O}OP}$ and by its coordinates
$x^1,\,x^2,\,x^3$:
$$
\hskip -2em
X^{i_1\ldots\,i_r}_{j_1\ldots\,j_s}=X^{i_1\ldots\,
i_r}_{j_1\ldots\,j_s}(x^1,x^2,x^3).
\tag19.3
$$
\proclaim{Conclusion 19.1} In contrast to free tensors, tensor
fields are related not to bases, but to whole coordinate systems
(including the origin). In each coordinate system they
are represented by functional arrays, i\.\,e\. by arrays of
functions (see \thetag{19.3}).
\endproclaim
    A functional array \thetag{19.3} is a coordinate representation
of a tensor field \thetag{19.1}. What happens when we change the
coordinate system\,? \pagebreak Dealing with \thetag{19.2}, we
need only to recalculate the components of the array $X^{i_1\ldots\,
i_r}_{j_1\ldots\,j_s}$ in the basis by means of \thetag{12.2}:
$$
\hskip -2em
\tilde X^{i_1\ldots\,i_r}_{j_1\ldots\,j_s}(P)=\msum\Sb h_1,\,\ldots,\,h_r\\
k_1,\,\ldots,\,k_s\endSb T^{i_1}_{h_1}\ldots\,T^{i_r}_{h_r}S^{k_1}_{j_1}
\ldots\,S^{k_s}_{j_s}\,X^{h_1\ldots\,h_r}_{k_1\ldots\,k_s}(P).
\tag19.4
$$
In the case of \thetag{19.3}, we need to recalculate the components of the
array $X^{i_1\ldots\,i_r}_{j_1\ldots\,j_s}$ in the new basis
$$
\tilde X^{i_1\ldots\,i_r}_{j_1\ldots\,j_s}(\tilde x^1,\tilde x^2,
\tilde x^3)=\msum\Sb h_1,\,\ldots,\,h_r\\
k_1,\,\ldots,\,k_s\endSb T^{i_1}_{h_1}\ldots\,T^{i_r}_{h_r}
S^{k_1}_{j_1}\ldots\,S^{k_s}_{j_s}\,X^{h_1\ldots\,h_r}_{k_1
\ldots\,k_s}(x^1,x^2,x^3),\ \quad
\tag19.5
$$
using \thetag{12.2}, and we also need to express the old coordinates
$x^1,\,x^2,\,x^3$ of the point $P$ in right hand side of \thetag{19.5}
through new coordinates of the same point:
$$
\hskip -2em
\cases
x^1=x^1(\tilde x^1,\tilde x^2,\tilde x^3),\\
x^2=x^2(\tilde x^1,\tilde x^2,\tilde x^3),\\
x^3=x^3(\tilde x^1,\tilde x^2,\tilde x^3).
\endcases
\tag19.6
$$
Like \thetag{12.2}, formula \thetag{19.5} can be inverted by
means of \thetag{12.1}:
$$
X^{i_1\ldots\,i_r}_{j_1\ldots\,j_s}(x_1,x_2,x_3)=
\msum\Sb h_1,\,\ldots,\,h_r\\
k_1,\,\ldots,\,k_s\endSb S^{i_1}_{h_1}\ldots\,
S^{i_r}_{h_r}T^{k_1}_{j_1}\ldots\,T^{k_s}_{j_s}
\,\tilde X^{h_1\ldots\,h_r}_{k_1\ldots\,k_s}(\tilde x^1,
\tilde x^2,\tilde x^3).\quad
\tag19.7
$$
But now, apart from \thetag{19.7}, we should have inverse formulas
for \thetag{19.6} as well:
$$
\hskip -2em
\cases
\tilde x^1=x^1(x^1,x^2,x^3),\\
\tilde x^2=x^2(x^1,x^2,x^3),\\
\tilde x^3=x^3(x^1,x^2,x^3).
\endcases
\tag19.8
$$
THe couple of formulas \thetag{19.5} and \thetag{19.6}, and another
couple of formulas \thetag{19.7} and \thetag{19.8}, in the case
of tensor fields play the same role as transformation formulas
\thetag{12.1} and \thetag{12.2} in the case of free tensors.
\head
\S~20. Change of Cartesian coordinate system.
\endhead
    Note that formulas \thetag{19.6} and \thetag{19.8} are written
in abstract form. They only indicate the functional dependence of
new coordinates of the point $P$ from old ones and vice versa.
Now we shall specify them for the case when one Cartesian coordinate
system is changed to another Cartesian coordinate system. Remember
that each Cartesian coordinate system is determined by some basis
and some fixed point (the origin). We consider two Cartesian coordinate
systems. Let the origins of the first and second systems be at the points
$O$ and $\tilde O$, respectively. Denote by $\bold e_1,\,\bold e_2,\,
\bold e_3$ the basis of the first coordinate system, and by
$\tilde\bold e_1,\,\tilde\bold e_2,\,\tilde\bold e_3$ the basis
of the second coordinate system (see Fig\.~7 below).\par
\parshape 3 0pt 360pt 0pt 360pt 180pt 180pt
    Let $P$ be some point in the space for whose coordinates we
are going to derive the specializations of formulas \thetag{19.6} and
\thetag{19.8}. Denote by $\bold r_{\sssize P}$ and $\tilde\bold
r_{\sssize P}$ the radius-vectors of this point in our two coordinate
\vadjust{\vskip 235pt\hbox to 0pt{\kern 0pt
\includegraphics{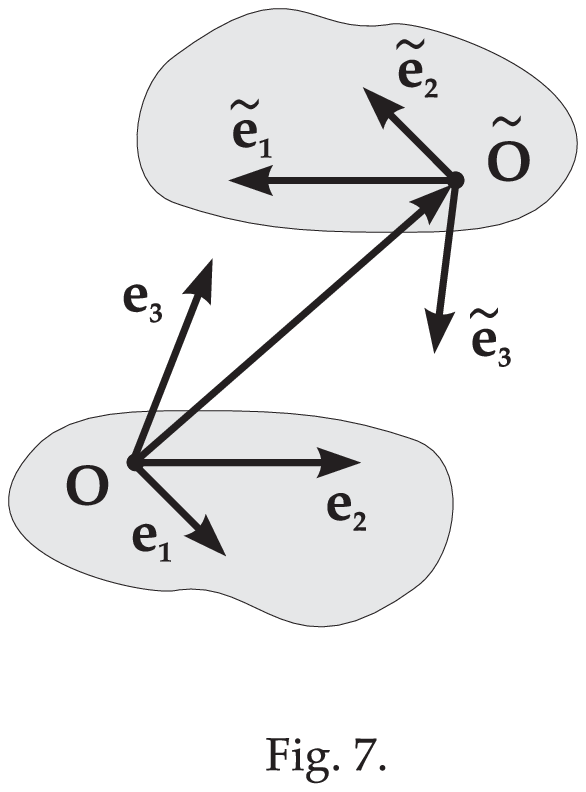}\hss}\vskip -235pt}systems.
Then $\bold r_{\sssize P}=
\overrightarrow{\vphantom{\tilde O}OP}$ and
$\tilde\bold r_{\sssize P}=\overrightarrow{\tilde OP}$. Hence
$$
\hskip -2em
\bold r_{\sssize P}=\overrightarrow{O\tilde O\vphantom{\vrule
height 11pt depth 1pt}}
+\tilde\bold r_{\sssize P}.
\tag20.1
$$
Vector $\overrightarrow{O\tilde O\vphantom{\vrule
height 11pt depth 1pt}}$ determines the origin shift from the old
to the new coordinate system. We expand this vector in the basis
$\bold e_1,\,\bold e_2,\,\bold e_3$:
$$
\hskip -2em
\bold a=\overrightarrow{O\tilde O\vphantom{\vrule height
11pt depth 1pt}}=\sum^3_{i=1}a^i\,\bold e_i.
\tag20.2
$$
Radius-vectors $\bold r_{\sssize P}$ and $\tilde\bold
r_{\sssize P}$ are expanded in the bases of their own
coordinate systems:
$$
\hskip -2em
\aligned
&\bold r_{\sssize P}=\sum^3_{i=1}x^i\,\bold e_i,\\
&\tilde\bold r_{\sssize P}=\sum^3_{i=1}\tilde x^i\,\tilde\bold e_i,
\endaligned
\tag20.3
$$
\proclaim{Exercise 20.1} Using \thetag{20.1}, \thetag{20.2},
\thetag{20.3}, and \thetag{5.7} derive the following formula
relating the coordinates of the point $P$ in the two coordinate
systems in Fig\.~7:
$$
\hskip -2em
x^i=a^i+\sum^3_{j=1}S^i_j\,\tilde x^j.
\tag20.4
$$
Compare \thetag{20.4} with \thetag{6.5}. Explain the differences
in these two formulas.
\endproclaim
\proclaim{Exercise 20.2} Derive the following inverse formula
for \thetag{20.4}:
$$
\hskip -2em
\tilde x^i=\tilde a^i+\sum^3_{j=1}T^i_j\,x^j.
\tag20.5
$$
Prove that $a^i$ in \thetag{20.4} and $\tilde a^i$ in \thetag{20.5}
are related to each other as follows:
$$
\xalignat 2
&\hskip -2em
\tilde a^i=-\sum^3_{j=1}T^i_j\,a^j,
&&a^i=-\sum^3_{j=1}S^i_j\,\tilde a^j.
\tag20.6
\endxalignat
$$
Compare \thetag{20.6} with \thetag{6.2} and \thetag{6.5}.
Explain the minus signs in these formulas.
\endproclaim
\noindent
Formula \thetag{20.4} can be written in the following
expanded form:
$$
\hskip -2em
\left\{\aligned
&x^1=S^1_1\,\tilde x^1+S^1_2\,\tilde x^2+S^1_3\,\tilde x^3+a^1,\\
&x^2=S^2_1\,\tilde x^1+S^2_2\,\tilde x^2+S^2_3\,\tilde x^3+a^2,\\
&x^3=S^3_1\,\tilde x^1+S^3_2\,\tilde x^2+S^3_3\,\tilde x^3+a^3.
\endaligned\right.
\tag20.7
$$
This is the required specialization for \thetag{19.6}.
In a similar way we can expand \thetag{20.5}:
$$
\hskip -2em
\left\{\aligned
&\tilde x^1=T^1_1\,x^1+T^1_2\,x^2+T^1_3\,x^3+\tilde a^1,\\
&\tilde x^2=T^2_1\,x^1+T^2_2\,x^2+T^2_3\,x^3+\tilde a^2,\\
&\tilde x^3=T^3_1\,x^1+T^3_2\,x^2+T^3_3\,x^3+\tilde a^3.
\endaligned\right.
\tag20.8
$$
This is the required specialization for \thetag{19.8}. Formulas
\thetag{20.7} and \thetag{20.8} are used to accompany the
main transformation formulas \thetag{19.5} and \thetag{19.7}.
\head
\S~21. Differentiation of tensor fields.
\endhead
    In this section we consider two different types of derivatives
that are usually applied to tensor fields: differentiation with
respect to spacial variables $x^1,\,x^2,\,x^3$ and differentiation
with respect to external parameters other than $x^1,\,x^2,\,x^3$,
if they are present. The second type of derivatives are simpler to
understand. Let's consider them to start. Suppose we have
tensor field $\bold X$ of type $(r,s)$ and depending on the
additional parameter $t$ (for instance, this could be a time variable).
Then, upon choosing some Cartesian coordinate system, we can write
$$
\frac{\partial X^{i_1\ldots\,i_r}_{j_1\ldots\,j_s}}{\partial t}
=\lim_{h\to 0}\frac{X^{i_1\ldots\,i_r}_{j_1\ldots\,j_s}(t+h,x^1,
x^2,x^3)-X^{i_1\ldots\,i_r}_{j_1\ldots\,j_s}(t,x^1,x^2,x^3)}{h}.
\quad
\tag21.1
$$
The left hand side of \thetag{21.1} is a tensor since the fraction
in right hand side is constructed by means of tensorial operations
\thetag{14.1} and \thetag{14.3}. Passing to the limit $h\to 0$ does
not destroy the tensorial nature of this fraction since the
transition matrices $S$ and $T$ in \thetag{19.5}, \thetag{19.7},
\thetag{20.7}, \thetag{20.8} are all time-independent.
\proclaim{Conclusion 21.1} Differentiation with respect to external
parameters (like $t$ in \thetag{21.1}) is a tensorial operation
producing new tensors from existing ones.
\endproclaim
\proclaim{Exercise 21.1} Give a more detailed explanation of why
the time derivative \thetag{21.1} represents a tensor of type
$(r,s)$.
\endproclaim
Now let's consider the spacial derivative of tensor field $\bold X$,
i\.\,e\. its derivative with respect to a spacial variable, e\.\,g.
with respect to $x^1$. Here we also can write
$$
\frac{\partial X^{i_1\ldots\,i_r}_{j_1\ldots\,j_s}}{\partial x^1}
=\lim_{h\to 0}\frac{X^{i_1\ldots\,i_r}_{j_1\ldots\,j_s}(x^1+h,
x^2,x^3)-X^{i_1\ldots\,i_r}_{j_1\ldots\,j_s}(x^1,x^2,x^3)}{h},
\quad
\tag21.2
$$
but in numerator of the fraction in the right hand side of
\thetag{21.2} we get the difference of two tensors bound to
different points of space: to the point $P$ with coordinates
$x^1,\,x^2,\,x^3$ and to the point $P'$ with coordinates
$x^1+h,\,x^2,\,x^3$. To which point should be attributed
the difference of two such tensors\,? This is not clear.
Therefore, we should treat partial derivatives like
\thetag{21.2} in a different way.\par
    Let's choose some additional symbol, say it can be $q$,
and consider the partial derivative of $X^{i_1\ldots\,i_r}_{j_1\ldots
\,j_s}$ with respect to the spacial variable $x^q$:
$$
Y^{i_1\ldots\,i_r}_{q\,j_1\ldots\,j_s}=
\frac{\partial X^{i_1\ldots\,i_r}_{j_1\ldots\,j_s}}{\partial x^q}.
\tag21.3
$$
Partial derivatives \thetag{21.2}, \pagebreak taken as a whole, form
an $(r+s+1)$-dimensional array with one extra dimension due
to index $q$. We write it as a lower index in $Y^{i_1\ldots\,
i_r}_{q\,j_1\ldots\,j_s}$ due to the following theorem
concerning \thetag{21.3}.
\proclaim{Theorem 21.1} For any tensor field $\bold X$ of type
$(r,s)$ partial derivatives \thetag{21.3} with respect to
spacial variables $x^1,\,x^2,\,x^3$ in any Cartesian coordinate
system represent another tensor field $\bold Y$ of the type
$(r,s+1)$.
\endproclaim
Thus differentiation with respect to $x^1,\,x^2,\,x^3$ produces
new tensors from already existing ones. For the sake of beauty and
convenience this operation is denoted by the nabla sign: $\bold Y=
\nabla\bold X$. In index form this looks like
$$
\hskip -2em
Y^{i_1\ldots\,i_r}_{q\,j_1\ldots\,j_s}=
\nabla_{\!q}X^{i_1\ldots\,i_r}_{j_1\ldots\,j_s}.
\tag21.4
$$
Simplifying the notations we also write:
$$
\hskip -2em
\nabla_{\!q}=\frac{\partial}{\partial x^q}.
\tag21.5
$$
\proclaim{Warning 21.1} Theorem~21.1 and the equality \thetag{21.5}
are valid only for Cartesian coordinate systems. In curvilinear
coordinates things are different.
\endproclaim
\proclaim{Exercise 21.2} Prove theorem~21.1. For this purpose
consider another Cartesian coordinate system $\tilde x^1,\,
\tilde x^2,\,\tilde x^3$ related to $x^1,x^2,x^3$ via
\thetag{20.7} and \thetag{20.8}. Then in the new coordinate system
consider the partial derivatives
$$
\hskip -2em
\tilde Y^{i_1\ldots\,i_r}_{q\,j_1\ldots\,j_s}=
\frac{\partial\tilde X^{i_1\ldots\,i_r}_{j_1\ldots\,j_s}}
{\partial\tilde x^q}
\tag21.6
$$
and derive relationships binding \thetag{21.6} and \thetag{21.3}.
\endproclaim
\head
\S~22. Gradient, divergency, and rotor.\\
Laplace and d'Alambert operators.
\endhead
\rightheadtext{\S~22. Gradient, divergency, and rotor \dots}
    The tensorial nature of partial derivatives established by
theorem~21.1 is a very useful feature. We can apply it to extend
the scope of classical operations of vector analysis.
Let's consider {\bf the gradient}, $\grad=\nabla$. Usually
the gradient operator is applied to scalar fields, i.\,e\. to
functions $\varphi=\varphi(P)$ or $\varphi=\varphi(x^1,x^2,
x^3)$ in coordinate form:
$$
\hskip -2em
a_q=\nabla_{\!q}\varphi=\frac{\partial\varphi}{\partial x^q}.
\tag22.1
$$
Note that in \thetag{22.1} we used a lower index $q$ for $a_q$.
This means that $\bold a=\grad\varphi$ is a covector. Indeed,
according to theorem~21.1, the nabla operator applied to a
scalar field, which is tensor field of type $(0,0)$, produces
a tensor field of type $(0,1)$. In order to get the vector form
of the gradient one should raise index $q$:
$$
\hskip -2em
a^q=\sum^3_{i=1}g^{qi}\,a_i=\sum^3_{i=1}g^{qi}\,\nabla_{\!i}
\varphi.
\tag22.2
$$
Let's write \thetag{22.2} in the form of a differential operator
(without applying to $\varphi$):
$$
\hskip -2em
\nabla^q=\sum^3_{i=1}g^{qi}\,\nabla_{\!i}.
\tag22.3
$$
In this form the gradient operator \thetag{22.3} can be applied
not only to scalar fields, but also to vector fields, covector
fields and to any other tensor fields.\par
    Usually in physics we do not distinguish between the
vectorial gradient $\nabla^q$ and the covectorial gradient
$\nabla_{\!q}$ because we use orthonormal coordinates with ONB as
a basis. In this case dual metric tensor is given by unit matrix
($g^{ij}=\delta^{ij}$)
and components of $\nabla^q$ and $\nabla_{\!q}$ coincide by value.
\par
    {\bf Divergency} is the second differential operation of vector
analysis. Usually it is applied to a vector field and is given by
formula:
$$
\hskip -2em
\divr\bold X=\sum^3_{i=1}\nabla_{\!i}X^i.
\tag22.4
$$
As we see, \thetag{22.4} is produced by contraction (see section~16)
from tensor $\nabla_{\!q}X^i$. Therefore we can generalize formula
\thetag{22.4} and apply divergency operator to arbitrary tensor field
with at least one upper index:
$$
\hskip -2em
(\divr\bold X)^{\ldots\,\ldots\,\ldots}_{\ldots\,\ldots\,\ldots}
=\sum^3_{s=1}\nabla_{\!s}X^{\ldots\,.s.\,\ldots}_{\ldots\,\ldots\,
\ldots}.
\tag22.5
$$
{\bf The Laplace operator} is defined as the divergency applied
to a vectorial gradient of something, it is denoted by the triangle
sign: $\triangle=\divr\grad$. From \thetag{22.3} and \thetag{22.5}
for Laplace operator $\triangle$ we derive the following formula:
$$
\hskip -2em
\triangle=\sum^3_{i=1}\sum^3_{j=1}g^{ij}\,\nabla_{\!i}\,\nabla_{\!j}.
\tag22.6
$$
Denote by $\square$ the following differential operator:
$$
\square=\frac{1}{c^2}\frac{\partial^2}{\partial t^2}-
\triangle.
\tag22.7
$$
Operator \thetag{22.7} is called {\bf the d'Alambert operator}
or {\bf wave operator}. In general relativity upon introducing
the additional coordinate $x^0=c\,t$ one usually rewrites the
d'Alambert operator in a form quite similar to \thetag{22.6}
(see my book \cite{5}, it is free for download from 

\blue{http:/\negskp/samizdat.mines.edu/}).\par
    And finally, let's consider {\bf the rotor operator} or
{\bf curl operator} (the term ``rotor'' is derived from
``rotation'' so that ``rotor'' and ``curl'' have approximately
the same meaning). The rotor operator is usually applied to
a vector field and produces another vector field: $\bold Y=\rot
\bold X$. Here is the formula for the $r$-th coordinate of
$\rot\bold X$:
$$
\hskip -2em
(\rot\bold X)^r=\sum^3_{i=1}\sum^3_{j=1}\sum^3_{k=1}g^{ri}
\,\omega_{ijk}\,\nabla^jX^k.
\tag22.8
$$
The volume tensor $\boldsymbol\omega$ in \thetag{22.8} is given
by formula \thetag{18.4}, while the vectorial gradient operator
$\nabla^j$ is defined in \thetag{22.3}.
\proclaim{Exercise 22.1} Formula \thetag{22.8} can be generalized
for the case when $\bold X$ is an arbitrary tensor field with at
least one upper index. By analogy with \thetag{22.5} suggest your
version of such a generalization.
\endproclaim
Note that formulas \thetag{22.6} and \thetag{22.8} for the Laplace
operator and for the rotor are different from those that are commonly
used. Here are standard formulas:
$$
\gather
\hskip -2em
\triangle=\left(\kern -1pt\frac{\partial}{\partial x^1}
\!\right)^{\kern-1pt\lower 2pt\hbox{$\ssize 2$}}+
\left(\kern -1pt\frac{\partial}{\partial x^2}\!
\right)^{\kern-1pt\lower 2pt\hbox{$\ssize 2$}}+
\left(\kern -1pt\frac{\partial}{\partial x^3}\!
\right)^{\kern-1pt\lower 2pt\hbox{$\ssize 2$}},
\tag22.9\\
\vspace{2ex}
\hskip -2em
\rot\bold X=
\det\Vmatrix \bold e_1 & \bold e_2 & \bold e_3\\
\vspace{2ex}
\dfrac{\partial}{\partial x^1} & \dfrac{\partial}
{\partial x^2} & \dfrac{\partial}{\partial x^3}\\
\vspace{2ex}
X^1 & X^2 & X^3
\endVmatrix.
\tag22.10
\endgather
$$
The truth is that formulas \thetag{22.6} and \thetag{22.8}
are written for a general skew-angular coordinate system with
a SAB as a basis. The standard formulas \thetag{22.10} are valid
only for orthonormal coordinates with ONB as a basis.
\proclaim{Exercise 22.2} Show that in case of orthonormal
coordinates, when $g^{ij}=\delta^{ij}$, formula \thetag{22.6}
for the Laplace operator $\triangle$ reduces to the standard
formula \thetag{22.9}.
\endproclaim
    The coordinates of the vector $\rot\bold X$ in a skew-angular
coordinate system are given by formula \thetag{22.8}. Then for
vector $\rot\bold X$ itself we have the expansion:
$$
\hskip -2em
\rot\bold X=\sum^3_{r=1}(\rot\bold X)^r\,\bold e_r.
\tag22.11
$$
\proclaim{Exercise 22.3} Substitute \thetag{22.8} into
\thetag{22.11} and show that in the case of a orthonormal
coordinate system the resulting formula \thetag{22.11}
reduces to \thetag{22.10}.
\endproclaim
\newpage
%--------------------------------------
\setfirstpage
\topmatter
\title\chapter{4}
TENSOR FIELDS IN CURVILINEAR COORDINATES.
\endtitle
\endtopmatter
\document
\head
\S~23. General idea of curvilinear coordinates.
\endhead
\leftheadtext{CHAPTER~\uppercase\expandafter{\romannumeral 4}.
TENSOR FIELDS IN CURVILINEAR COORDINATES.}
    What are coordinates, if we forget for a moment about
radius-vectors, bases and axes\,? What is the pure idea of
coordinates\,? The pure idea is in representing points of space
by triples of numbers. This means that we should have one to
one map $P\leftrightarrows (y^1,y^2,y^3)$ in the whole space
or at least in some domain, where we are going to use our
coordinates $y^1,\,y^2,\,y^3$. In Cartesian coordinates this
map $P\leftrightarrows (y^1,y^2,y^3)$ is constructed by means
of vectors and bases. Arranging other coordinate systems one
can use other methods. For example, in {\bf spherical coordinates}
$y^1=r$ is a distance from the point $P$ to the center of
sphere, $y^2=\theta$ and $y^3=\varphi$ are two angles. By
the way, spherical coordinates are one of the simplest
examples of curvilinear coordinates. Furthermore, let's keep in
mind spherical coordinates when thinking about more general
and hence more abstract curvilinear coordinate systems.
\head
\S~24. Auxiliary Cartesian coordinate system.
\endhead
    Now we know almost everything about Cartesian coordinates
and almost nothing about the abstract curvilinear coordinate system
$y^1,\,y^2,\,y^3$ that we are going to study. Therefore, the best
idea is to represent each point $P$ by its radius-vector
$\bold r_{\!\sssize P}$ in some auxiliary Cartesian coordinate
system and then consider a map $\bold r_{\!\sssize P}\leftrightarrows
(y^1,y^2,y^3)$. The radius-vector itself is represented by three
coordinates in the basis $\bold e_1,\,\bold e_2,\,\bold e_3$
of the auxiliary Cartesian coordinate system:
$$
\hskip -2em
\bold r_{\!\sssize P}=\sum^3_{i=1}x^i\,\bold e_i.
\tag24.1
$$
Therefore, we have a one-to-one map $(x^1,x^2,x^3)\leftrightarrows
(y^1,y^2,y^3)$. Hurrah! This is a numeric map. We can treat it
numerically. In the left direction it is represented by three
functions of three variables:
$$
\hskip -2em
\cases
x^1=x^1(y^1,y^2,y^3),\\
x^2=x^2(y^1,y^2,y^3),\\
x^3=x^3(y^1,y^2,y^3).
\endcases
\tag24.2
$$
In the right direction we again have three functions of three variables:
$$
\hskip -2em
\cases
y^1=y^1(x^1,x^2,x^3),\\
y^2=y^2(x^1,x^2,x^3),\\
y^3=y^3(x^1,x^2,x^3).
\endcases
\tag24.3
$$
Further we shall assume all functions in \thetag{24.2} and
\thetag{24.3} to be differentiable and consider their partial
derivatives. Let's denote
$$
\xalignat 2
&\hskip -2em
S^i_j=\frac{\partial x^i}{\partial y^j},
&&T^i_j=\frac{\partial y^i}{\partial x^j}.
\tag24.4
\endxalignat
$$
Partial derivatives \thetag{24.4} can be arranged into
two square matrices $S$ and $T$ respectively. In mathematics
such matrices are called Jacobi matrices. The components of
matrix $S$ in that form, as they are defined in \thetag{24.4},
are functions of $y^1,\,y^2,\,y^3$. The components of $T$ are
functions of $x^1,\,x^2,\,x^3$:
$$
\xalignat 2
&\hskip -2em
S^i_j=S^i_j(y^1,y^2,y^3),
&&T^i_j=T^i_j(x^1,x^2,x^3).
\tag24.5
\endxalignat
$$
However, by substituting \thetag{24.3} into the arguments of
$S^i_j$, or by substituting \thetag{24.2} into the arguments
of $T^i_j$, we can make them have a common set of arguments:
$$
\xalignat 2
&\hskip -2em
S^i_j=S^i_j(x^1,x^2,x^3),
&&T^i_j=T^i_j(x^1,x^2,x^3),
\tag24.6\\
&\hskip -2em
S^i_j=S^i_j(y^1,y^2,y^3),
&&T^i_j=T^i_j(y^1,y^2,y^3),
\tag24.7
\endxalignat
$$
When brought to the form \thetag{24.6}, or when brought to
the form \thetag{24.7} (but not in form of \thetag{24.5}),
matrices $S$ and $T$ are inverse of each other: 
$$
\hskip -2em
T=S^{-1}.
\tag24.8
$$
This relationship \thetag{24.8} is due to the fact that
numeric maps \thetag{24.2} and \thetag{24.3} are inverse
of each other.\par
\proclaim{Exercise 24.1} You certainly know the following
formula:
$$
\gather
\frac{df(x^1(y),x^2(y),x^3(y))}{dy}=
\sum^3_{i=1}f'_i(x^1(y),x^2(y),x^3(y))
\,\frac{dx^i(y)}{dy},
\text{\ \ where \ }f'_i=\frac{\partial f}{\partial x^i}.
\endgather
$$
It's for the differentiation of composite function. Apply this formula
to functions \thetag{24.2} and derive the relationship \thetag{24.8}.
\endproclaim
\head
\S~25. Coordinate lines and the coordinate grid.
\endhead
    Let's substitute \thetag{24.2} into \thetag{24.1} and take
into account that \thetag{24.2} now assumed to contain differentiable
functions. Then the vector-function 
$$
\hskip -2em
\bold R(y^1,y^2,y^3)=\bold r_{\!\sssize P}=\sum^3_{i=1}
x^i(y^1,y^2,y^3)\,\bold e_i
\tag25.1
$$
is a differentiable function of three variables $y^1,\,y^2,\,y^3$.
The vector-function $\bold R(y^1,y^2,y^3)$ determined by \thetag{25.1}
is called {\bf a basic vector-function} of a curvilinear coordinate
system. Let $P_0$ be some fixed point of space given by its
curvilinear coordinates $y^1_0,\,y^2_0,\,y^3_0$. Here zero is not the
tensorial index, we use it in order to emphasize that $P_0$ is fixed
point, and that its coordinates $y^1_0,\,y^2_0,\,y^3_0$ are
three fixed numbers. In the next step let's undo one of them,
say first one, by setting
$$
\xalignat 3
&\hskip -2em
y^1=y^1_0+t, &&y^2=y^2_0, &&y^3=y^3_0.
\tag25.2
\endxalignat
$$\par
\parshape 19 0pt 360pt 180pt 180pt 180pt 180pt 180pt 180pt
180pt 180pt 180pt 180pt 180pt 180pt 180pt 180pt 180pt 180pt
180pt 180pt 180pt 180pt 180pt 180pt 180pt 180pt 180pt 180pt
180pt 180pt 180pt 180pt 180pt 180pt 180pt 180pt
0pt 360pt
\noindent
Substituting \thetag{25.2} into \thetag{25.1} we get
a vector-function \vadjust{\vskip 210pt\hbox to 0pt{\kern 0pt
\includegraphics{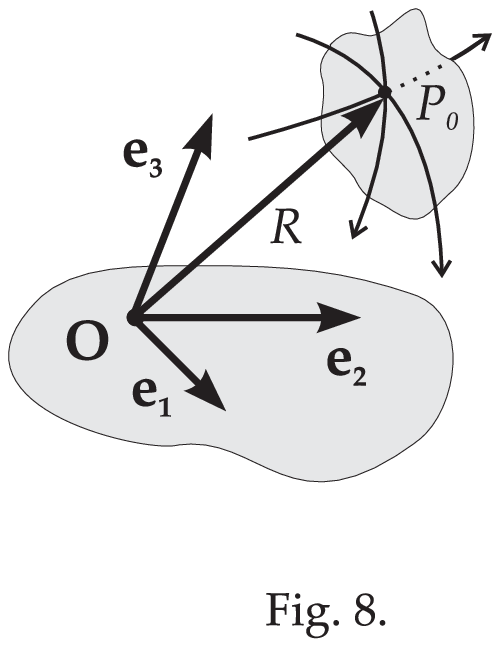}\hss}\vskip -210pt}of
one variable $t$:
$$
\bold R_1(t)=\bold R(y^1_0+t,y^2_0,y^3_0),\quad
\tag25.3
$$
If we treat $t$ as time variable (though it may have a
unit other than time), then \thetag{25.3} describes a curve
(the trajectory of a particle). At time instant $t=0$ this
curve passes through the fixed point $P_0$. Same is true for
curves given by two other vector-functions similar to
\thetag{25.4}:
$$
\align
&\hskip -4em
\bold R_2(t)=\bold R(y^1_0,\,y^2_0+t,\,y^3_0),
\tag25.4\\
\vspace{1ex}
&\hskip -4em
\bold R_3(t)=\bold R(y^1_0,\,y^2_0,\,y^3_0+t).
\tag25.5
\endalign
$$
This means that all three curves given by
vector-functions \thetag{25.3}, \thetag{25.4}, and
\thetag{25.5} are intersected
at the point $P_0$ as shown on Fig\.~8.
\vadjust{\vskip 190pt\hbox to 0pt{\kern 190pt\includegraphics{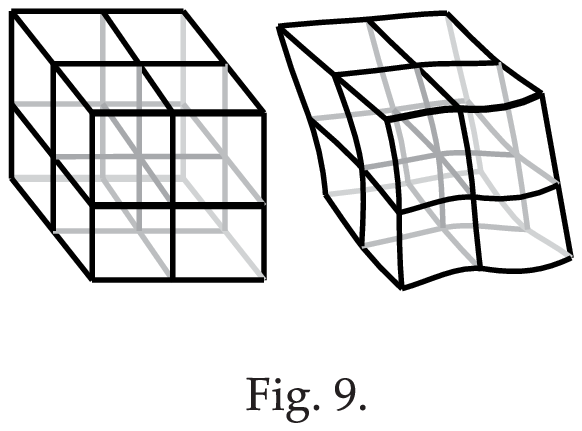}\hss}\vskip -190pt}Arrowheads on these lines
indicate the directions in which parameter $t$ increases. Curves
\thetag{25.3}, \thetag{25.4}, and \thetag{25.5} are called
{\bf coordinate lines}. They are subdivided into three families.
Curves within one family do not intersect each other. Curves from
different families intersect so that any regular point of
space is an intersection of exactly three coordinate curves (one
per family).\par
\parshape 1 0pt 180pt
    Coordinate lines taken in whole form {\bf a coordinate grid}.
This is an infinitely dense grid. But usually, when drawing, it is
represented as a grid with finite density. On Fig\.~9 the coordinate
grid of curvilinear coordinates is compared to that of the Cartesian
coordinate system.\par
\parshape 6 0pt 180pt 0pt 180pt 0pt 180pt 0pt 180pt 0pt 180pt
0pt 360pt
    Another example of coordinate grid is on Fig\.~2. Indeed, meridians
and parallels are coordinate lines of a spherical coordinate system.
The parallels do not intersect, but the meridians
forming one family of coordinate lines do intersect at the North
and at South Poles. This means that North and South Poles are
singular points for spherical coordinates.
\proclaim{Exercise 25.1} Remember the exact definition of spherical
coordinates and find all singular points for them. 
\endproclaim
\head
\S~26. Moving frame of curvilinear coordinates.
\endhead
    Let's consider the three coordinate lines shown on Fig.~8 again.
And let's find tangent vectors to them at the point $P_0$. For
this purpose we should differentiate vector-functions
\thetag{25.3}, \thetag{25.4}, and \thetag{25.5} with respect
to the time variable $t$ and then substitute $t=0$ into the
derivatives:
$$
\bold E_i=\frac{d\bold R_i}{dt}\,\hbox{\vrule height 12pt depth 8pt
width 0.5pt}_{\,t=0}=\frac{\partial\bold R}{\partial y^i}\hbox{\vrule
height 12pt depth 8pt width 0.5pt}_{\,\text{at the point $P_0$}}.
\tag26.1
$$
Now let's substitute the expansion \thetag{25.1} into \thetag{26.1}
and remember \thetag{24.4}:
$$
\bold E_i=\frac{\partial\bold R}{\partial y^i}=
\sum^3_{j=1}\frac{\partial x^j}{\partial y^i}\,\bold e_j=
\sum^3_{j=1}S^j_i\,\bold e_j.
\tag26.2
$$
All calculations in \thetag{26.2} are still in reference to the
point $P_0$. Though $P_0$ is a fixed point, it is an arbitrary
fixed point. Therefore, the equality \thetag{26.2} is valid at
any point. Now let's omit the intermediate calculations and write
\thetag{26.2} as
$$
\hskip -2em
\bold E_i=\sum^3_{i=1}S^j_i\,\bold e_j.
\tag26.3
$$\par
\parshape 3 0pt 360pt 0pt 360pt 180pt 180pt 
\noindent
And then compare \thetag{26.3} with \thetag{5.7}. They are strikingly
similar, and $\det S\neq 0$\linebreak due to \thetag{24.8}. Formula
\thetag{26.3}
\vadjust{\vskip 210pt\hbox to 0pt{\kern 0pt\includegraphics{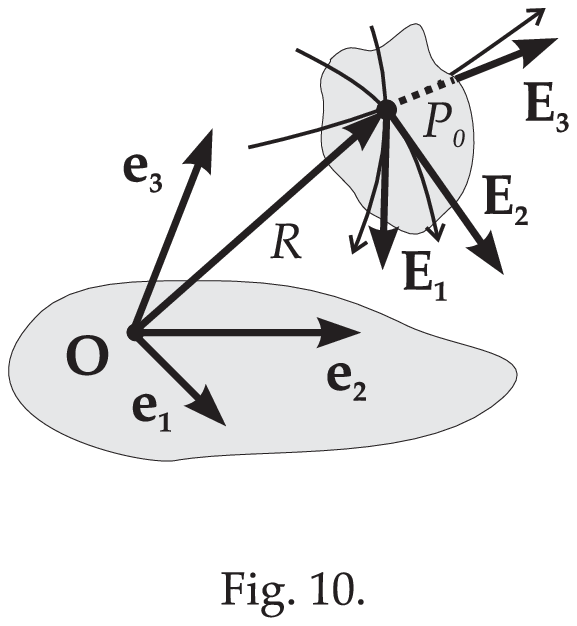}\hss}\vskip -210pt}means that tangent vectors to coordinate
lines $\bold E_1,\,\bold E_2,\,\bold E_3$ form a basis (see Fig\.~10),
matrices \thetag{24.4} are transition matrices to this basis and
back to the Cartesian basis.\par
\parshape 14 180pt 180pt 180pt 180pt 180pt 180pt 180pt 180pt
180pt 180pt 180pt 180pt 180pt 180pt 180pt 180pt 180pt 180pt
180pt 180pt 180pt 180pt 180pt 180pt 180pt 180pt 0pt 360pt
    Despite obvious similarity of the formulas \thetag{26.3} and
\thetag{5.7}, there is some crucial difference of basis $\bold E_1,
\,\bold E_2,\,\bold E_3$ as compa\-red to $\bold e_1,\bold e_2,\,
\bold e_3$. Vectors $\bold E_1,\,\bold E_2,$ $\bold E_3$ are not
free. They are bound to that point where derivatives \thetag{24.4}
are calculated. And they move when we move this point. For this
reason basis $\bold E_1,\,\bold E_2,\,\bold E_3$ is called
{\bf moving frame} of the curvilinear coordinate system. During
their motion the vectors of the moving frame $\bold E_1,\,
\bold E_2,\,\bold E_3$ are not simply translated\linebreak from
point to
point, they can change their lengths and the angles they form
with each other. Therefore, in general the moving frame
$\bold E_1,\,\bold E_2,\,\bold E_3$ is a skew-angular basis.
In some cases vectors $\bold E_1,\,\bold E_2,\,\bold E_3$
can be orthogonal to each other at all points of space. In that case
we say that we have an orthogonal curvilinear coordinate system.
Most of the well known curvilinear coordinate systems are
orthogonal, e\.\,g\. spherical, cylindrical, elliptic, parabolic,
toroidal, and others. However, there is no curvilinear coordinate
system with the moving frame being ONB\,! We shall not prove this
fact since it leads deep into differential geometry.
\head
\S~27. Dynamics of moving frame.
\endhead
   Thus, we know that {\bf the moving frame} moves. Let's
describe this motion quantitatively. According to \thetag{24.5}
the components of matrix $S$ in \thetag{26.3} are functions of
the curvilinear coordinates $y^1,\,y^2,\,y^3$. Therefore,
differentiating $\bold E_i$ with respect to $y^j$ we should
expect to get some nonzero vector $\partial\bold E_i/
\partial y^j$. This vector can be expanded back in moving frame
$\bold E_1,\,\bold E_2,\,\bold E_3$. This expansion is written as
$$
\hskip -2em
\frac{\partial\bold E_i}{\partial y^j}=\sum^3_{k=1}
\Gamma^k_{ij}\,\bold E_k.
\tag27.1
$$
Formula \thetag{27.1} is known as {\bf the derivational formula}.
Coefficients $\Gamma^k_{ij}$ in \thetag{27.1} are called {\bf
Christoffel symbols} or {\bf connection components}.\par
\proclaim{Exercise 27.1} Relying upon formula \thetag{25.1} and
\thetag{26.1} draw the vectors of the moving frame for cylindrical
coordinates.
\endproclaim
\proclaim{Exercise 27.2} Do the same for spherical coordinates.
\endproclaim
\proclaim{Exercise 27.3} Relying upon formula \thetag{27.1} and
results  of exercise~27.1 calculate the Christoffel symbols for
cylindrical coordinates.
\endproclaim
\proclaim{Exercise 27.4} Do the same for spherical coordinates.
\endproclaim
\proclaim{Exercise 27.5} Remember formula \thetag{26.2} from which
you derive
$$
\hskip -2em
\bold E_i=\frac{\partial\bold R}{\partial y^i}.
\tag27.2
$$
Substitute \thetag{27.2} into left hand side of the derivational
formula \thetag{27.1} and relying on the properties of mixed
derivatives prove that the Christoffel symbols are symmetric
with respect to their lower indices:
$\Gamma^k_{ij}=\Gamma^k_{ji}$.
\endproclaim
    Note that Christoffel symbols $\Gamma^k_{ij}$ form a
three-dimensional array with one upper index and two lower
indices. However, they do not represent a tensor. We shall
not prove this fact since it again leads deep into differential
geometry.
\head
\S~28. Formula for Christoffel symbols.
\endhead
    Let's take formula \thetag{26.3} and substitute it into
both sides of \thetag{27.1}. As a result we get the following
equality for Christoffel symbols $\Gamma^k_{ij}$:
$$
\hskip -2em
\sum^3_{q=1}\frac{\partial S^q_i}{\partial y^j}\,\bold e_q=
\sum^3_{k=1}\sum^3_{q=1}\Gamma^k_{ij}\,S^q_k\,\bold e_q.
\tag28.1
$$
Cartesian basis vectors $\bold e_q$ do not depend on $y^j$;
therefore, they are not differentiated when we substitute
\thetag{26.3} into \thetag{27.1}. Both sides of \thetag{28.1}
are expansions in the base $\bold e_1,\,\bold e_2,\,\bold e_3$
of the auxiliary Cartesian coordinate system. Due to the
uniqueness of such expansions we have the following equality
derived from \thetag{28.1}:
$$
\hskip -2em
\frac{\partial S^q_i}{\partial y^j}=
\sum^3_{k=1}\Gamma^k_{ij}\,S^q_k.
\tag28.2
$$
\proclaim{Exercise 28.1} Using concept of the inverse matrix
(\,$T=S^{-1}$) derive the following formula for the Christoffel
symbols $\Gamma^k_{ij}$ from \thetag{28.2}:
$$
\hskip -2em
\Gamma^k_{ij}=
\sum^3_{q=1}T^k_q\,\frac{\partial S^q_i}{\partial y^j}.
\tag28.3
$$
\endproclaim
Due to \thetag{24.4} this formula \thetag{28.3} can be transformed
in the following way:
$$
\hskip -2em
\Gamma^k_{ij}=
\sum^3_{q=1}T^k_q\,\frac{\partial S^q_i}{\partial y^j}=
\sum^3_{q=1}T^k_q\,\frac{\partial^2 x^q}{\partial y^i\,
\partial y^j}=\sum^3_{q=1}T^k_q\,\frac{\partial S^q_j}
{\partial y^i}.
\tag28.4
$$
Formulas \thetag{28.4} are of no practical use because
they express $\Gamma^k_{ij}$ through an external thing like
transition matrices to and from the auxiliary Cartesian
coordinate system. However, they will help us below in
understanding the differentiation of tensors.
\head
\S~29. Tensor fields in curvilinear coordinates. 
\endhead
    As we remember, tensors are geometric objects related
to bases and represented by arrays if some basis is specified.
Each curvilinear coordinate system provides us a numeric
representation for points, and in addition to this it provides
the basis. This is the moving frame. Therefore, we can refer
tensorial objects to curvilinear coordinate systems, where they
are represented as arrays of functions:
$$
\hskip -2em
X^{i_1\ldots\,i_r}_{j_1\ldots\,j_s}=
X^{i_1\ldots\,i_r}_{j_1\ldots\,j_s}(y^1,y^2,y^3).
\tag29.1
$$
We also can have two curvilinear coordinate systems and can
pass from one to another by means of transition functions:
$$
\xalignat 2
&\hskip -2em
\cases
\tilde y^1=\tilde y^1(y^1,y^2,y^3),\\
\tilde y^2=\tilde y^2(y^1,y^2,y^3),\\
\tilde y^3=\tilde y^3(y^1,y^2,y^3),
\endcases
&&\cases
y^1=y^1(\tilde y^1,\tilde y^2,\tilde y^3),\\
y^2=y^2(\tilde y^1,\tilde y^2,\tilde y^3),\\
y^3=y^3(\tilde y^1,\tilde y^2,\tilde y^3).
\endcases
\tag29.2
\endxalignat
$$
If we call $\tilde y^1,\,\tilde y^2,\,\tilde y^3$ the new coordinates,
and $y^1,\,y^2,\,y^3$ the old coordinates, then transition matrices
$S$ and $T$ are given by the following formulas:
$$
\xalignat 2
&\hskip -2em
S^i_j=\frac{\partial y^i}{\partial\tilde y^j},
&&T^i_j=\frac{\partial\tilde y^i}{\partial y^j}.
\tag29.3
\endxalignat
$$
They relate moving frames of two curvilinear coordinate
systems:
$$
\xalignat 2
&\hskip -2em
\tilde\bold E_i=\sum^3_{j=1}S^{\,j}_i\,\bold E_j,
&&\bold E_j=\sum^3_{i=1}T^i_j\,\tilde \bold E_i.
\tag29.4
\endxalignat
$$
\proclaim{Exercise 29.1} Derive \thetag{29.3} from
\thetag{29.4} and \thetag{29.2} using some auxiliary
Cartesian coordinates with basis $\bold e_1,\,\bold e_2,
\bold e_3$ as intermediate coordinate system:
\vskip -2.4ex
$$
\hskip -2em
(\bold E_1,\,\bold E_2,\,\bold E_3)
\vcenter{\hsize 45pt$$\CD @<S<<\\ \vspace{-1.5em} @>>T>\endCD$$}
(\bold e_1,\,\bold e_2,\,\bold e_3)
\vcenter{\hsize 45pt$$\CD @>\tilde S>>\\ \vspace{-1.5em}
@<<\tilde T<\endCD$$} (\tilde\bold E_1,\,\tilde\bold E_2,\,
\tilde\bold E_3)
\tag29.5
$$
\vskip -2ex
\noindent
Compare \thetag{29.5} with \thetag{5.13} and explain differences
you have detected.
\endproclaim
Transformation formulas for tensor fields for two curvilinear 
coordinate systems are the same as in \thetag{19.4} and \thetag{19.5}:
$$
\gather
\hskip -3em
\tilde X^{i_1\ldots\,i_r}_{j_1\ldots\,j_s}(\tilde y^1,\tilde y^2,
\tilde y^3)=\msum\Sb h_1,\,\ldots,\,h_r\\
k_1,\,\ldots,\,k_s\endSb T^{i_1}_{h_1}\ldots\,T^{i_r}_{h_r}
S^{k_1}_{j_1}\ldots\,S^{k_s}_{j_s}\,X^{h_1\ldots\,h_r}_{k_1
\ldots\,k_s}(y^1,y^2,y^3),\ \quad
\tag29.6\\
\hskip -3em
X^{i_1\ldots\,i_r}_{j_1\ldots\,j_s}(y_1,y_2,y_3)=
\msum\Sb h_1,\,\ldots,\,h_r\\
k_1,\,\ldots,\,k_s\endSb S^{i_1}_{h_1}\ldots\,
S^{i_r}_{h_r}T^{k_1}_{j_1}\ldots\,T^{k_s}_{j_s}
\,\tilde X^{h_1\ldots\,h_r}_{k_1\ldots\,k_s}(\tilde y^1,
\tilde y^2,\tilde y^3).\quad
\tag29.7
\endgather
$$
But formulas \thetag{19.6} and \thetag{19.8} should be
replaced by \thetag{29.2}.
\head
\S~30. Differentiation of tensor fields\\
in curvilinear coordinates.
\endhead
    We already know how to differentiate tensor fields
in Cartesian coordinates (see section~21). We know that
operator $\nabla$ produces tensor field of type $(r,s+1)$
when applied to a tensor field of type $(r,s)$. The only
thing we need now is to transform $\nabla$ to a curvilinear
coordinate system. In order to calculate tensor $\nabla\bold X$
in curvilinear coordinates, let's first transform $\bold X$ into
auxiliary Cartesian coordinates, then apply $\nabla$, and
then transform $\nabla\bold X$ back into curvilinear coordinates:
$$
\hskip -2em
\CD
X^{h_1\ldots\,h_r}_{k_1\ldots\,k_s}(y^1,y^2,y^3)
@>S,T>>X^{h_1\ldots\,h_r}_{k_1\ldots\,k_s}(x^1,x^2,x^3)\\
@VV\nabla_{\!p} V @VV\nabla_{\!q}=\partial/\partial x^q V\\
\nabla_{\!p}X^{i_1\ldots\,i_r}_{j_1\ldots\,j_s}(y^1,y^2,y^3)
@<T,S<<\nabla_{\!q}X^{h_1\ldots\,h_r}_{k_1\ldots\,k_s}(x^1,x^2,x^3)
\endCD
\tag30.1
$$
Matrices \thetag{24.4} are used in \thetag{30.1}. From \thetag{12.3}
and \thetag{12.4} we know that the transformation of each index is
a separate multiplicative procedure. When applied to the $\alpha$-th
upper index, the whole chain of transformations \thetag{30.1} looks
like
$$
\nabla_{\!p}X^{\ldots\,i_{\!\sssize\alpha}\,\ldots}_{\ldots\,
\ldots\,\ldots}=\sum^3_{q=1}S^q_p\ldots\sum^3_{h_{\!\sssize\alpha}=1}
T^{\,i_{\!\sssize\alpha}}_{h_{\!\sssize\alpha}}\,\ldots\,\nabla_{\!q}\,
\ldots\sum^3_{m_{\!\sssize\alpha}=1}S^{\,h_{\!\sssize\alpha}}_{m_{\!
\sssize\alpha}}\ldots\, X^{\ldots\,m_{\!\sssize\alpha}\,
\ldots}_{\ldots\,\ldots\,\ldots}.\quad
\tag30.2
$$
Note that $\nabla_{\!q}=\partial/\partial x^q$ is a differential operator
and due to \thetag{24.4} we have
$$
\hskip -2em
\sum^3_{q=1}S^q_p\,\frac{\partial}{\partial x^q}=
\frac{\partial}{\partial y^p}.
\tag30.3
$$
Any differential operator when applied to a product produces a sum
with as many summands as there were multiplicand in the product.
Here is the summand produced by term $S^{\,h_{\!\sssize\alpha}}_{m_{\!
\sssize\alpha}}$ in formula \thetag{30.2}:
$$
\hskip -2em
\nabla_{\!p}X^{\ldots\,i_{\!\sssize\alpha}\,\ldots}_{\ldots\,
\ldots\,\ldots}=\ldots+\sum^3_{m_{\!\sssize\alpha}=1}
\sum^3_{h_{\!\sssize\alpha}=1}T^{\,i_{\!\sssize\alpha}}_{h_{\!
\sssize\alpha}}\,\frac{S^{\,h_{\!\sssize\alpha}}_{m_{\!\sssize
\alpha}}}{\partial y^p}\,X^{\ldots\,m_{\!\sssize\alpha}\,
\dots}_{\ldots\,\ldots\,\ldots}+\ldots\,.
\tag30.4
$$
Comparing \thetag{30.4} with \thetag{28.3} or \thetag{28.4}
we can transform it into the following equality:
$$
\hskip -2em
\nabla_{\!p}X^{\ldots\,i_{\!\sssize\alpha}\,\ldots}_{\ldots\,
\ldots\,\ldots}=\ldots+\sum^3_{m_{\!\sssize\alpha}=1}
\Gamma^{\,i_{\!\sssize\alpha}}_{pm_{\!
\sssize\alpha}}\,X^{\ldots\,m_{\!\sssize\alpha}\,\dots}_{\ldots
\,\ldots\,\ldots}+\ldots\,.
\tag30.5
$$
Now let's consider the transformation of the $\alpha$-th lower
index in \thetag{30.1}:
$$
\nabla_{\!p}X^{\ldots\,\ldots\,\ldots}_{\ldots\,j_{\!\sssize\alpha}\,
\ldots}=\sum^3_{q=1}S^q_p\ldots\sum^3_{k_{\!\sssize\alpha}=1}
S^{k_{\!\sssize\alpha}}_{j_{\!\sssize\alpha}}\,\ldots\,\nabla_{\!q}\,
\ldots\sum^3_{n_{\!\sssize\alpha}=1}T^{n_{\!\sssize\alpha}}_{k_{\!
\sssize\alpha}}\ldots\,X^{\ldots\,\ldots\,\ldots}_{\ldots\,n_{\!
\sssize\alpha}\,\ldots}.\quad
\tag30.6
$$
Applying \thetag{30.3} to \thetag{30.6} with the same logic as
in deriving \thetag{30.4} we get
$$
\hskip -2em
\nabla_{\!p}X^{\ldots\,\ldots\,\ldots}_{\ldots\,j_{\!\sssize\alpha}\,
\ldots}=\ldots+\sum^3_{n_{\!\sssize\alpha}=1}
\sum^3_{k_{\!\sssize\alpha}=1}S^{k_{\!\sssize\alpha}}_{j_{\!\sssize
\alpha}}\,\frac{T^{n_{\!\sssize\alpha}}_{k_{\!\sssize\alpha}}}
{\partial y^p}\,X^{\ldots\,\ldots\,\ldots}_{\ldots\,n_{\!
\sssize\alpha}\,\ldots}+\ldots\,.
\tag30.7
$$
In order to simplify \thetag{30.7} we need the following formula
derived from \thetag{28.3}:
$$
\hskip -2em
\Gamma^k_{ij}=
-\sum^3_{q=1}S^q_i\,\frac{\partial T^k_q}{\partial y^j}.
\tag30.8
$$
Applying \thetag{30.8} to \thetag{30.7} we obtain
$$
\hskip -2em
\nabla_{\!p}X^{\ldots\,\ldots\,\ldots}_{\ldots\,j_{\!\sssize\alpha}\,
\ldots}=\ldots-\sum^3_{n_{\!\sssize\alpha}=1}
\Gamma^{n_{\!\sssize\alpha}}_{pj_{\!\sssize\alpha}}\,X^{\ldots\,
\ldots\,\ldots}_{\ldots\,n_{\!\sssize\alpha}\,\ldots}+\ldots\,.
\tag30.9
$$
Now we should gather \thetag{30.5}, \thetag{30.9}, and 
add the term produced when $\nabla_{\!q}$ in \thetag{30.2}
(or equivalently in \thetag{30.4}) acts upon components of
tensor $\bold X$. As a result we get the following general
formula for $\nabla_{\!p}X^{i_1\ldots\,i_r}_{j_1\ldots\,j_s}$:
$$
\hskip -2em
\gathered
\nabla_{\!p}X^{i_1\ldots\,i_r}_{j_1\ldots\,j_s}=
\frac{\partial X^{i_1\ldots\,i_r}_{j_1\ldots\,j_s}}{\partial y^p}+
\sum^r_{\alpha=1}\sum^3_{m_{\!\sssize\alpha}=1}
\Gamma^{\,i_{\!\sssize\alpha}}_{pm_{\!\sssize\alpha}}
\,X^{i_1\ldots\,m_{\!\sssize\alpha}\,\ldots\,i_r}_{
j_1\ldots\,\,\ldots\,\,\ldots\, j_s}-\\
-\sum^s_{\alpha=1}\sum^3_{n_{\!\sssize\alpha}=1}
\Gamma^{n_{\!\sssize\alpha}}_{pj_{\!\sssize\alpha}}
\,X^{i_1\ldots\,\,\ldots\,\,\ldots\,i_r}_{j_1\ldots
\,n_{\!\sssize\alpha}\,\ldots\,j_s}.
\endgathered
\tag30.10
$$
The operator $\nabla_p$ determined by this formula is called
{\bf the covariant derivative}. 
\proclaim{Exercise 30.1} Apply the general formula \thetag{30.10}
to a vector field and calculate the covariant derivative $\nabla_{\!p}X^q$.
\endproclaim
\proclaim{Exercise 30.2} Apply the general formula \thetag{30.10}
to a covector field and calculate the covariant derivative
$\nabla_{\!p}X_q$.
\endproclaim
\proclaim{Exercise 30.3} Apply the general formula \thetag{30.10}
to an operator field and find $\nabla_{\!p}F^q_m$. Consider
special case when $\nabla_{\!p}$ is applied to the Kronecker symbol
$\delta^q_m$.
\endproclaim
\proclaim{Exercise 30.4} Apply the general formula \thetag{30.10}
to a bilinear form and find $\nabla_{\!p}a_{qm}$.
\endproclaim
\proclaim{Exercise 30.5} Apply the general formula \thetag{30.10}
to a tensor product $\bold a\otimes\bold x$ for the case when
$\bold x$ is a vector and $\bold a$ is a covector. Verify formula
$\nabla(\bold a\otimes\bold x)=\nabla\bold a\otimes\bold x+
\bold a\otimes\nabla\bold x$.
\endproclaim
\proclaim{Exercise 30.6} Apply the general formula
\thetag{30.10} to the contraction $C(\bold F)$ for
the case when $\bold F$ is an operator field.
Verify the formula $\nabla C(\bold F)=C(\nabla\bold F)$.
\endproclaim
\proclaim{Exercise 30.7} Derive \thetag{30.8} from \thetag{28.3}.
\endproclaim
\head
\S~31. Concordance of metric and connection.
\endhead
    Let's remember that we consider curvilinear coordinates in
Euclidean space $E$. In this space we have the scalar product
\thetag{13.1} and the metric tensor \thetag{13.5}.
\proclaim{Exercise 31.1} Transform the metric tensor \thetag{13.5}
to curvilinear coordinates using transition matrices \thetag{24.4}
and show that here it is given by formula
$$
\hskip -2em
g_{ij}=(\bold E_i,\,\bold E_j).
\tag31.1
$$
\endproclaim
In Cartesian coordinates all components of the metric tensor are constant
since the basis vectors $\bold e_1,\,\bold e_2,\,\bold e_3$ are constant.
The covariant derivative \thetag{30.10} in Cartesian coordinates reduces
to differentiation $\nabla_{\!p}=\partial/\partial x^p$. Therefore,
$$
\hskip -2em
\nabla_{\!p}g_{ij}=0.
\tag31.2
$$
But $\nabla g$ is a tensor. If all of its components in some coordinate
system are zero, then they are identically zero in any other coordinate
system (explain why). Therefore the identity \thetag{31.2} is valid in
curvilinear coordinates as well.
\proclaim{Exercise 31.2} Prove \thetag{31.2} by direct calculations
using formula \thetag{27.1}.
\endproclaim
\noindent
The identity \thetag{31.2} is known as {\bf the concordance condition}
for the metric $g_{ij}$ and connection $\Gamma^k_{ij}$. It is
very important for general relativity.\par
    Remember that the metric tensor enters into many useful formulas
for the gradient, divergency, rotor, and Laplace operator in
section 22. What is important is that all of these formulas remain
valid in curvilinear coordinates, with the only difference being
that you should understand that $\nabla_{\!p}$ is not the partial
derivative $\partial/\partial x^p$, but the covariant derivative in
the sense of formula \thetag{30.10}.
\proclaim{Exercise 31.3} Calculate $\rot\bold A$, $\divr\bold H$,
$\grad\varphi$ (vectorial gradient) in cylindrical and spherical
coordinates.
\endproclaim
\proclaim{Exercise 31.4} Calculate the Laplace operator $\triangle
\varphi$ applied to the scalar field $\varphi$ in cylindrical and
in spherical coordinates.
\endproclaim
\newpage
%--------------------------------------
\topmatter
\title
REFERENCES.
\endtitle
\endtopmatter
\document
\setfirstpage

\enddocument
\end